\documentclass{amsart}
\usepackage{amssymb}

\newcommand{\newpar}{}

\newcommand{\CalP}{\mathcal{P}}
\newcommand{\CalH}{\mathcal{H}}

\newcommand{\Z}{\mathbb{Z}}

\newcommand{\T}{\mathcal{T}}
\newcommand{\Fq}{\mathbb{F}_q}
\newcommand{\Fqs}{\mathbb{F}_{q^{2}}}
\newcommand{\Fqe}{\mathbb{F}_{q^{e}}}
\newcommand{\Qlbar}{\overline{\mathbb{Q}_l}}

\newcommand{\Ind}{\mathrm{Ind}}
\newcommand{\tr}{\mathrm{tr}}
\newcommand{\inv}{\mathrm{inv}}
\newcommand{\ff}{\mathrm{ff}}
\newcommand{\tsigma}{\tilde{\sigma}}
\newcommand{\isomto}{\overset{\sim}{\rightarrow}}

\numberwithin{equation}{subsection}

\theoremstyle{plain}
\newtheorem{theorem}{Theorem}[subsection]
\newtheorem{corollary}[theorem]{Corollary}
\newtheorem{lemma}[theorem]{Lemma}
\newtheorem{proposition}[theorem]{Proposition}

\theoremstyle{definition}

\begin{document}
\title[$K^{F}$-invariants in irreducible representations
of $G^{F}$, when $G=GL_n$]
{$K^{F}$-invariants in irreducible \\ representations
of $G^{F}$, when $G=GL_n$}

\author{Anthony Henderson}
\address{School of Mathematics and Statistics,
University of Sydney, NSW 2006, AUSTRALIA}
\email{anthonyh@maths.usyd.edu.au}
\begin{abstract}
Using a general result of Lusztig, we give explicit
formulas for the dimensions of $K^{F}$-invariants in
irreducible representations of $G^{F}$, when $G=GL_n$,
$F:G\to G$ is a Frobenius map, and $K$ is an $F$-stable
subgroup of finite index in $G^{\theta}$ for some involution
$\theta:G\to G$ commuting with $F$. The proofs use
some combinatorial facts about characters of symmetric groups.
\end{abstract}
\maketitle
\section*{Introduction}
\newpar
Let $G$ be a connected reductive group defined over a finite field
$\Fq$ of odd characteristic. Let $F$ be the Frobenius morphism on $G$,
whose fixed-point subgroup $G^{F}=G(\Fq)$ is finite.
Let $\theta:G\to G$ be an involution of algebraic groups
commuting with $F$, and $K$ an $F$-stable subgroup of
the fixed-point subgroup $G^{\theta}$ which contains
$(G^{\theta})^{\circ}$. The homogeneous variety
$G/K$ is a \emph{symmetric space}, and the set of cosets
$G^{F}/K^{F}$ might reasonably be called a \emph{finite symmetric space}.

\newpar
The classification of irreducible representations of $G^{F}$
(in characteristic $0$) was completed by Lusztig in the mid-1980s
(see \cite{lusztigwarsaw} for a general statement, and
\cite{chars} and \cite{disconnected} for more details). A roughly analogous
problem for symmetric spaces is that of determining the dimension
of the $K^{F}$-invariants in each irreducible representation of $G^{F}$;
equivalently, calculating the multiplicities
$\langle\chi, \Ind_{K^{F}}^{G^{F}}(1)\rangle$
for every irreducible character $\chi$ of $G^{F}$. A large first
step towards solving this problem was Lusztig's calculation, in
\cite{symmfinite}, of
$\langle\tr(\cdot,R_T^{\lambda}), \Ind_{K^{F}}^{G^{F}}(1)\rangle$
for every Deligne-Lusztig virtual representation $R_T^{\lambda}$ of $G^{F}$.
In \cite{gfq-inv}, Lusztig proceeded to
solve the problem completely in the case when
$G^{F}=(G')^{F^{2}}$, $K^{F}=(G')^{F}$, and $G'$ has connected centre.

\newpar
This paper is devoted to the solution of this problem when $G$
is a general linear group (with either split or non-split $F$,
so that $G^{F}$ is either $GL_n(\Fq)$ or $U_n(\Fqs)$)
and $\theta$ is arbitrary. (The solution for the case $G^{F}=GL_{n}(\Fq)$,
$K^{F}=Sp_{n}(\Fq)$ was found by Bannai, Kawanaka, and Song
in \cite[\S4]{bks}.)
For such $G$, the functions $\tr(\cdot,R_T^{\lambda})$ form a basis of the
class functions, and the transition matrix from this basis to that of the
irreducible characters is known. So Lusztig's result gives
a formula for $\langle\chi, \Ind_{K^{F}}^{G^{F}}(1)\rangle$. 
All that remains is to manipulate
this formula until it is manifestly a nonnegative integer,
a straightforward (though not entirely easy) matter. 
Two justifications for presenting
it in detail are the potential interest of the answers, and the pleasantness
of the symmetric group combinatorics involved.

\newpar
In \S1 we recall Lusztig's formula and the character theory of the
finite general linear and unitary groups as well as introducing
some vital notation. Then we traverse
the various cases in \S\!\S2-4, which could be thought of as a theme
and variations: the theme, or underlying pattern, is stated in its
simplest form in \S2.1 (the case already known from \cite{bks}),
and successive subsections follow the same pattern with progressively
more elaborate alterations. 
In \S2 the involution is symplectic, so $n$ must be even and
the possible symmetric spaces are $GL_n(\Fq)/Sp_n(\Fq)$ and
$U_n(\Fqs)/Sp_n(\Fq)$. In \S3 the involution is inner,
so the possible symmetric spaces are $GL_n(\Fq)/
(GL_{n^{+}}(\Fq)\times GL_{n^{-}}(\Fq))$,
$GL_n(\Fq)/GL_{n/2}(\Fqs)$, $U_n(\Fqs)/
(U_{n^{+}}(\Fqs)\times U_{n^{-}}(\Fqs))$, and
$U_n(\Fqs)/U_{n/2}(\mathbb{F}_{q^{4}})$. In \S4 the involution is
orthogonal, so $G^{\theta}$ is not connected. However as noted
in Lemma \ref{sonlemma}, it is enough to solve the problem
when $K=G^{\theta}$, in which case the possible symmetric spaces are
$GL_n(\Fq)/O_n(\Fq)$ ($n$ odd), $GL_n(\Fq)/O_n^{\pm}(\Fq)$ ($n$ even),
$U_n(\Fqs)/O_n(\Fq)$ ($n$ odd), and $U_n(\Fqs)/O_n^{\pm}(\Fq)$ ($n$ even).

\newpar
The key combinatorial results we need along
the way are all proved in \S5. A reader interested only in these results
could skip all of \S\!\S1-4 except \S1.2.

\newpar
To give some idea
of how the formulas in \S\!\S2-4 connect with previously known results, 
we here extract the
answers for \emph{unipotent} irreducible characters. For both
$GL_n(\Fq)$ and $U_n(\Fqs)$, these are parametrized by partitions of $n$,
say $\rho\mapsto\chi^{\rho}\in\widehat{G^{F}}$. (In our convention
$\chi^{(n)}$ is the trivial character and $\chi^{(1^{n})}$ is the
Steinberg character.) Recall that a signed tableau of shape
$\mu$ is a signed Young diagram of shape $\mu$ where signs alternate
across rows, modulo permutations of rows of equal length. For this and 
all other combinatorial notation, see \S1.2. We have:
\[ (n \text{ even})\qquad
\langle \chi^{\rho}, \Ind_{Sp_n(\Fq)}^{GL_n(\Fq)}(1)\rangle=
\left\{
\begin{array}{cl}
1, &\text{ if $\rho$ is even}\\ 
0, &\text{ otherwise} 
\end{array} \right. \]
(for the general $GL_n(\Fq)/Sp_n(\Fq)$ case,
see Theorem \ref{glnspnthm});
\[ (n \text{ even})\qquad
\langle \chi^{\rho}, \Ind_{Sp_n(\Fq)}^{U_n(\Fqs)}(1)\rangle=
\left\{
\begin{array}{cl}
1, &\text{ if $\rho$ is even}\\ 
0, &\text{ otherwise} 
\end{array} \right. \]
(for the general $U_n(\Fqs)/Sp_n(\Fq)$ case,
see Theorem \ref{unspnthm});
\[ \langle \chi^{\rho},
\Ind_{GL_{n^{+}}(\Fq)\times GL_{n^{-}}(\Fq)}^{GL_n(\Fq)}(1)\rangle=
\begin{array}{l}
\text{the number of signed tableaux of shape $\rho'$}\\
\text{and signature $(n^{+},n^{-})$}
\end{array} \]
(for the general $GL_n(\Fq)/
(GL_{n^{+}}(\Fq)\times GL_{n^{-}}(\Fq))$ case,
see Theorem \ref{glnglnglnthm});
\begin{equation*}
\begin{split}
(n \text{ even})\qquad
\langle \chi^{\rho},
\Ind_{GL_{n/2}(\Fqs)}^{GL_n(\Fq)}(1)\rangle&=
\begin{array}{l}
\text{the number of signed tableaux of shape $\rho'$,}\\
\text{stable under changing all signs}
\end{array} \\
&=\left\{
\begin{array}{cl}
1, &\text{ if $\rho$ is even}\\ 
0, &\text{ otherwise} 
\end{array} \right.
\end{split}
\end{equation*}
(for the general $GL_n(\Fq)/
GL_{n/2}(\Fqs)$ case,
see Theorem \ref{glnglnthm});
\[ \langle \chi^{\rho},
\Ind_{U_{n^{+}}(\Fqs)\times U_{n^{-}}(\Fqs)}^{U_n(\Fqs)}(1)\rangle=
\begin{array}{l}
\text{the number of signed tableaux of shape $\rho'$}\\
\text{and signature $(n^{+},n^{-})$, stable under}\\
\text{inverting all rows}
\end{array} \]
(for the general $U_n(\Fqs)/
(U_{n^{+}}(\Fqs)\times U_{n^{-}}(\Fqs))$ case,
see Theorem \ref{unununthm});
\begin{equation*}
\begin{split}
(n \text{ even})\qquad
\langle \chi^{\rho},
\Ind_{U_{n/2}(\mathbb{F}_{q^{4}})}^{U_n(\Fqs)}(1)\rangle&=
\begin{array}{l}
\text{the number of signed tableaux of shape $\rho'$,}\\
\text{stable under changing all signs}\\
\text{and inverting all rows}
\end{array} \\
&=\left\{
\begin{array}{cl}
{\displaystyle \prod_i (m_{2i}(\rho')+1)}, 
&\text{ if $2\,|\,m_{2i+1}(\rho'),\ \forall i$}\\ 
0, &\text{ otherwise}
\end{array} \right.
\end{split}
\end{equation*}
(for the general $U_n(\Fqs)/
U_{n/2}(\mathbb{F}_{q^{4}})$ case,
see Theorem \ref{ununthm});
\begin{equation*}
(n \text{ odd})\qquad
\langle \chi^{\rho},
\Ind_{O_n(\Fq)}^{GL_n(\Fq)}(1)\rangle
=\frac{1}{2}\prod_i (m_{i}(\rho)+1)
\end{equation*}
(for the general $GL_n(\Fq)/
O_n(\Fq)$ case, see Theorem \ref{glnonthm});
\begin{equation*}
(n \text{ even})\qquad
\langle \chi^{\rho},
\Ind_{O_n^{\epsilon}(\Fq)}^{GL_n(\Fq)}(1)\rangle=
\left\{\begin{array}{cl}
{\displaystyle\lceil\frac{1}{2}\prod_i (m_{i}(\rho)+1)\rceil}, 
&\text{ if $\rho'$ is even, $\epsilon=+$}\\
{\displaystyle\lfloor\frac{1}{2}\prod_i (m_{i}(\rho)+1)\rfloor}, 
&\text{ if $\rho'$ is even, $\epsilon=-$}\\ 
{\displaystyle\frac{1}{2}\prod_i (m_{i}(\rho)+1)}, &\text{ otherwise} 
\end{array} \right.
\end{equation*}
(for the general $GL_n(\Fq)/
O_n^{\pm}(\Fq)$ case, see Theorem \ref{glnonpmthm});
\begin{equation*}
(n \text{ odd})\qquad
\langle \chi^{\rho},
\Ind_{O_n(\Fq)}^{U_n(\Fqs)}(1)\rangle
=\left\{\begin{array}{cl}
{\displaystyle\frac{1}{2}\prod_i (m_{2i+1}(\rho)+1)}, 
&\text{ if $2\,|\,m_{2i}(\rho),\ \forall i$}\\
0, &\text{ otherwise}
\end{array}\right.
\end{equation*}
(for the general $U_n(\Fqs)/
O_n(\Fq)$ case, see Theorem \ref{unonthm}); and
\begin{equation*}
(n \text{ even})\qquad
\langle \chi^{\rho},
\Ind_{O_n^{\epsilon}(\Fq)}^{U_n(\Fqs)}(1)\rangle=
\left\{\begin{array}{cl}
{\displaystyle\lceil\frac{1}{2}\prod_i (m_{2i}(\rho)+1)\rceil},
&\text{ if $\rho'$ is even, $\epsilon=+$}\\
{\displaystyle\lfloor\frac{1}{2}\prod_i (m_{2i}(\rho)+1)\rfloor},
&\text{ if $\rho'$ is even, $\epsilon=-$}\\
{\displaystyle\frac{1}{2}\prod_i (m_{2i}(\rho)+1)},
&\begin{array}{l}
\text{if $2\,|\,m_{2i+1}(\rho),\ \forall i$,}\\
\text{but $\rho'$ not even}\end{array}\\
0, &\text{ otherwise}
\end{array}\right.
\end{equation*}
(for the general $U_n(\Fqs)/
O_n^{\pm}(\Fq)$ case, see Theorem \ref{unonpmthm}).

\newpar
Readers experienced in the theory of cells for the symmetric group
will find these answers familiar. In fact, most of 
the above facts about unipotent
characters can be obtained by a more direct method than the one used
in this paper. For instance, suppose that $G^{F}\cong GL_n(\Fq)$.
Then the unipotent irreducible characters are the constituents of
$\Ind_{B^{F}}^{G^{F}}(1)$ where $B$ is an $F$-stable Borel subgroup.
The Hecke algebra $\CalH(G^{F},B^{F})$
is the specialization at $q$ of the abstract Hecke algebra $\CalH$ of $S_n$;
let $V_\rho$ be the simple $\CalH$-module indexed by $\rho$.

\newpar
It is trivial to show that the above multiplicity
$\langle \chi^{\rho}, \Ind_{K^{F}}^{G^{F}}(1)\rangle$
equals the multiplicity of $(V_\rho)_q$ in the $\CalH(G^{F},B^{F})$-module
$\mathcal{C}(B^{F}\!\setminus\! G^{F}/K^{F})$ of functions
on $G^{F}$ which are constant on the $B^{F}$--$K^{F}$ double cosets.
Assuming that $K$ is connected and split over $\Fq$, this module
is the specialization at $q$ of the $\CalH$-module $M^{K}$
defined in \cite{lusztigvogan}. So
$\langle \chi^{\rho}, \Ind_{K^{F}}^{G^{F}}(1)\rangle$
is the multiplicity of $V_\rho$ in $M^{K}$. In the current type-$A$ case,
this equals the number of \emph{cells} of $M^{K}$ which afford the
representation $V_\rho$, for which there is a combinatorial formula.
For example, our answer in the case of $GL_n(\Fq)/
(GL_{n^{+}}(\Fq)\times GL_{n^{-}}(\Fq))$ could be deduced from the
results in \cite{barbaschvogan}.
When $G$ and $K$ are non-split, this argument must be refined
to incorporate folding involutions, in the manner of \cite[\S10]{leftcells}.

\newpar
In principle, such arguments apply to
$\langle \chi^{\underline{\rho}}, \Ind_{K^{F}}^{G^{F}}(1)\rangle$
whenever $\chi^{\underline{\rho}}$ is a constituent of
$\mathrm{Ind}_{B^{F}}^{G^{F}}(\lambda)$,
since the generalized Hecke algebras 
$\mathrm{End}_{G^{F}}\mathrm{Ind}_{B^{F}}^{G^{F}}(\lambda)$
have been completely described. But the requisite facts about cells
for these Hecke algebras are somewhat diffuse in the
literature, and usually quoted in the slightly different context of real
Lie groups. I hope that the results of this paper, which are deduced
from \cite{symmfinite} in an independent way, will in fact shed further light
on the theory of cells.

\newpar
\textit{Acknowledgements.}
Part of this work was done at the Massachusetts Institute of Technology.
I would like to thank George Lusztig for all his help and encouragement,
and Peter Trapa and David Vogan for stimulating discussions.
I am also grateful to Dipendra Prasad and Kevin McGerty for impelling
me to revise the first version of this paper.
\tableofcontents
\section{Review of Known Results}
In \S1.1 we recall the theorem of Lusztig which underpins
all our results, and in \S1.3 the relevant parts of the
character theory of the finite general linear
and unitary groups in a convenient form. 
In \S1.2 and \S1.4, we introduce some combinatorial notation
to be used throughout the paper.
\subsection{Lusztig's Formula}
Let $k$ be the algebraic closure of a finite field $\Fq$ of odd cardinality
$q$. Let $G$ be a connected reductive group over $k$ defined over $\Fq$, with
Frobenius map $F:G\to G$. Let $\theta:G\to G$ be an involution of 
algebraic groups
commuting with $F$, and $K$ an $F$-stable subgroup of 
the fixed-point subgroup $G^{\theta}$ which contains
$(G^{\theta})^{\circ}$.
Fix a prime $l$ not dividing $q$. All representations and characters
of finite groups in this paper will be over $\Qlbar$.

\newpar
A \emph{pair} $(T,\lambda)$ means an $F$-stable maximal torus $T$
and a character $\lambda:T^{F}\to\Qlbar^{\times}$. 
We have a conjugation action of $G^{F}$ on the set of pairs.
In \cite{delignelusztig},
Deligne and Lusztig attached to each pair a virtual representation
$R_T^{\lambda}$ of $G^{F}$, depending only on the $G^{F}$-orbit of
$(T,\lambda)$. (In general, ``most'' of the $R_T^{\lambda}$ are irreducible
up to sign, and every irreducible representation occurs in some
$R_T^{\lambda}$; when $G=GL_n$ the situation is even better, as we will see
below.) The main result of \cite{symmfinite} is a formula for
\[ \frac{1}{|K^{F}|}\sum_{k\in K^{F}} \mathrm{tr}(k,R_T^{\lambda})
= \langle\tr(\cdot,R_T^{\lambda}), \Ind_{K^{F}}^{G^{F}}(1)\rangle. \]
After some trivial adjustments, it reads as follows. Define
\[ \Theta_T = \{f\in G\,|\, \theta(f^{-1}Tf)=f^{-1}Tf\}. \]
Then $T$ acts on $\Theta_T$ by left multiplication and $K$ acts
by right multiplication. 
If $B$ is a Borel subgroup containing $T$, the obvious map
$T\!\setminus\!\Theta_T/K\to B\!\setminus\!G/K$ is a bijection
(see \cite[Proposition 1.3]{symmfinite}), so $T\!\setminus\!\Theta_T/K$
is in bijection with the set of $K$-orbits on the flag variety.
For any $f\in\Theta_T^{F}$, define
$\epsilon_{T,f}:(T\cap fKf^{-1})^{F}\to\{\pm 1\}$ by
\[ \epsilon_{T,f}(t)=(-1)^{
\text{$\Fq$-rank}(Z_G((T\cap fKf^{-1})^{\circ}))
+\text{$\Fq$-rank}(Z_G^{\circ}(t)\cap Z_G((T\cap fKf^{-1})^{\circ}))}. \]
It follows from \cite[Proposition 2.3]{symmfinite} that
$\epsilon_{T,f}$ is a group homomorphism which factors through
$(T\cap fKf^{-1})^{F}/((T\cap fKf^{-1})^{\circ})^{F}$. Finally, define
\[ \Theta_{T,\lambda}^{F}=\{f\in\Theta_T^{F}\,|\,
\lambda|_{(T\cap fKf^{-1})^{F}}=\epsilon_{T,f}\}, \]
a union of $T^{F}$--$K^{F}$ double cosets.
\begin{theorem} \label{lusztigthm}
\textup{(Lusztig, \cite[Theorem 3.3]{symmfinite})}
\[ \langle\tr(\cdot,R_T^{\lambda}), \Ind_{K^{F}}^{G^{F}}(1)\rangle
=\sum_{f\in T^{F}\setminus \Theta_{T,\lambda}^{F}/K^{F}}
(-1)^{\textup{$\Fq$-rank}(T)+\textup{$\Fq$-rank}
(Z_G((T\cap fKf^{-1})^{\circ}))}.
\]
\end{theorem}
\subsection{Combinatorial Notation}
In general, our combinatorial notation always follows \cite{macdonald}.
For instance, $\mu\vdash n$ means that $\mu$ is a partition
of $n$. The size of a partition $\mu$ is written $|\mu|$ and its length
$\ell(\mu)$; it has parts $\mu_1$, $\mu_2$, \dots,
$\mu_{\ell(\mu)}$. The transpose partition is
$\mu'$.  We define
\[ n(\mu)=\sum_i (i-1)\mu_i =\sum_i\binom{\mu_i'}{2}. \]
The multiplicity of $i$ as a part of $\mu$ is written
$m_i(\mu)$. We say that $\mu$ is \emph{even} if all its parts
are even, or equivalently if $2|m_{i}(\mu'),\ \forall i$.

\newpar
It will be useful to have, for any partition $\nu$, a concrete
realization of the symmetric group $S_{|\nu|}$ which includes
a canonical element of cycle-type $\nu$. Let $\Lambda(\nu)$
be a set indexing the parts of $\nu$; usually,
$\Lambda(\nu)=\{1,\cdots,\ell(\nu)\}$. Then the set
\[ \Omega(\nu) = \{(j,i)\,|\,j\in\Lambda(\nu), i\in\Z/\nu_j\Z\} \]
has $|\nu|$ elements, and we will write $S_{|\nu|}$ for the group
of permutations of $\Omega(\nu)$ (a slight abuse of notation, since
it is not canonically associated to $|\nu|$). Let $w_\nu$ be the permutation
$(j,i)\mapsto (j,i+1)$. By construction this has cycle-type $\nu$,
and in fact we can identify $\Lambda(\nu)$ with the set of cycles of
$w_\nu$. Write $Z^{\nu}$ for the centralizer $Z_{S_{|\nu|}}(w_\nu)$.
The sign of $w_\nu$ is written $\epsilon_\nu$, and the size of
$Z^{\nu}$ is $z_\nu$. 

\newpar
For any $w\in Z^\nu$, let $\bar{w}$ be the induced permutation of
$\Lambda(\nu)$. Note that $\nu_{\bar{w}(j)}=\nu_{j}$ always.
For $w\in Z^\nu$ and $j\in\Lambda(\nu)$, there is a unique
$i(w,j)\in\Z/\nu_j\Z$ such that
\[ w(j,i)=(\bar{w}(j),i+i(w,j)),\ \forall i\in\Z/\nu_j\Z. \]
We will have much to do
with the set $Z^{\nu}_\inv$ of \emph{involutions} in $Z^{\nu}$.
Also define $Z^{\nu}_{\ff-\inv}$, the set of fixed-point free involutions.
Clearly $w\in Z^\nu$ lies in $Z^{\nu}_\inv$ if and only if $\bar{w}$
is an involution and
\[ i(w,\bar{w}(j))=-i(w,j),\ \forall j\in\Lambda(\nu). \]

\newpar
Note that any $w\in Z^{\nu}_\inv$ divides the set $\Lambda(\nu)$ of cycles
of $w_\nu$ into three disjoint subsets $\Lambda_w^{1}(\nu)$,
$\Lambda_w^{2}(\nu)$, and $\Lambda_w^{3}(\nu)$, according to whether
$w$ fixes the cycle pointwise, fixes the cycle but not pointwise,
or does not fix the cycle. Explicitly,
\begin{equation*}
\begin{split}
\Lambda_w^{1}(\nu)&=\{j\in\Lambda(\nu)\,|\,\bar{w}(j)=j, i(w,j)=0\},\\
\Lambda_w^{2}(\nu)&=\{j\in\Lambda(\nu)\,|\,\bar{w}(j)=j, i(w,j)=
\frac{1}{2}\nu_j\},\\
\Lambda_w^{3}(\nu)&=\{j\in\Lambda(\nu)\,|\,\bar{w}(j)\neq j\}.
\end{split}
\end{equation*}
We will write $\ell_w^{1}(\nu)$, $\ell_w^{2}(\nu)$, and
$\ell_w^{3}(\nu)$ for $|\Lambda_w^{1}(\nu)|$, $|\Lambda_w^{2}(\nu)|$, 
and $|\Lambda_w^{3}(\nu)|$, so that
\[ \ell_w^{1}(\nu) + \ell_w^{2}(\nu) + \ell_w^{3}(\nu) = \ell(\nu). \]
Note that $\ell_w^{3}(\nu)$ is always even, and
$w\in Z^{\nu}_{\ff-\inv}\Leftrightarrow\ell_w^{1}(\nu)=0$.

\newpar
We will use the notations
$\ell(\nu)_0$ and $\ell(\nu)_1$ for the number of even and odd
parts respectively. For instance, $\epsilon_\nu=(-1)^{\ell(\nu)_0}$.
On occasion we will need to further analyse $\ell(\nu)_0$
into $\ell(\nu)_{0\,\mathrm{mod}\,4}$ and $\ell(\nu)_{2\,\mathrm{mod}\,4}$.
We will also combine these notations in the obvious way,
e.g.\ $\ell_w^{3}(\nu)_1$ means the number of odd cycles of $w_\nu$
moved by $w$, and $\ell_w^{2}(\nu)_1=0$ always.

\newpar
In \S3 we will need to consider involutions with signed fixed points.
Let $(S_{|\nu|})_{\pm-\inv}$ be the set of pairs
$(w,\epsilon)$ where $w\in(S_{|\nu|})_{\inv}$ and
$\epsilon:\{(j,i)\,|\,w(j,i)=(j,i)\}\to\{+,-\}$ is a way
of signing the fixed points of $w$.
We define the \emph{signature} of $(w,\epsilon)$ to be
$(|\epsilon^{-1}(+)|,|\epsilon^{-1}(-)|)$.
We declare that
signatures are considered as elements of $\Z^{2}/\Z(1,1)$, e.g.\
$(1,0)$ and $(2,1)$ are the same.
Let $(S_{|\nu|})_{(p^{+},p^{-})-\inv}$ be the set of 
$(w,\epsilon)\in(S_{|\nu|})_{\pm-\inv}$ with signature $(p^{+},p^{-})$.
Now define
\[ Z^{\nu}_{\pm-\inv}=\{(w,\epsilon)\in(S_{|\nu|})_{\pm-\inv}\,|\,
w\in Z^{\nu}, \epsilon\circ w_\nu=\epsilon\}. \]
Note that
\[ |Z^{\nu}_{\pm-\inv}|=\sum_{w\in Z^{\nu}_\inv}2^{\ell_w^{1}(\nu)}. \]
Define $Z^{\nu}_{(p^{+},p^{-})-\inv}$ similarly, and also
\[ Z^{\nu}_{\star-\inv}=\{(w,\epsilon)\in(S_{|\nu|})_{\pm-\inv}\,|\,
w\in Z^{\nu}, \epsilon\circ w_\nu=-\epsilon\}. \]
Note that any $(w,\epsilon)\in Z^{\nu}_{\star-\inv}$ must have signature
$(0,0)$. Also
\[ |Z^{\nu}_{\star-\inv}|=\sum_{\substack{w\in Z^{\nu}_\inv\\
\ell_w^{1}(\nu)_1=0}}2^{\ell_w^{1}(\nu)}. \]

\newpar
Let $\T_{\pm}(\mu)$ be the set of \emph{signed tableaux of shape $\mu$}.
These are ways of labelling the boxes of the Young diagram of $\mu$
with a sign, in such a way that signs alternate across
each row (so all signs in a row are determined by that of the last box),
with the proviso that two labellings which differ by a permutation
of rows of equal length are not distinguished. Hence
\[ |\T_{\pm}(\mu)|=\prod_i (m_i(\mu)+1). \]
For $T\in\T_{\pm}(\mu)$, the \emph{signature} of $T$, again in
$\Z^{2}/\Z(1,1)$, is defined as
\[ (|\{\text{boxes of $T$ signed }+\}|,|\{\text{boxes of $T$ signed }-\}|), \]
or equivalently
\[ (|\{\text{odd rows of $T$ ending }\boxplus\}|,
|\{\text{odd rows of $T$ ending }\boxminus\}|). \]
Write $\T_{(p^{+},p^{-})}(\mu)$ for the set of $T\in\T_{\pm}(\mu)$
with signature $(p^{+},p^{-})$. There are two important involutions on
$\T_{\pm}(\mu)$: $\varphi$ which changes all signs, and
$\psi$ which reverses all rows. Write $\T_{\pm}(\mu)^\varphi$ etc.\
for the fixed-point sets. Note that any $T$ in
$\T_{\pm}(\mu)^\varphi$ or $\T_{\pm}(\mu)^{\varphi\psi}$ must have
signature $(0,0)$. Also
\begin{equation*}
\begin{split}
|\T_{\pm}(\mu)^\varphi|
&=\left\{
\begin{array}{cl}
1, &\text{ if $2\,|\,m_{i}(\mu),\ \forall i$}\\ 
0, &\text{ otherwise,}
\end{array} \right.\\
|\T_{\pm}(\mu)^\psi|
&=\left\{
\begin{array}{cl}
{\displaystyle \prod_i (m_{2i+1}(\mu)+1)}, 
&\text{ if $2\,|\,m_{2i}(\rho'),\ \forall i$}\\ 
0, &\text{ otherwise, and}
\end{array} \right. \\
|\T_{\pm}(\mu)^{\varphi\psi}|
&=\left\{
\begin{array}{cl}
{\displaystyle \prod_i (m_{2i}(\mu)+1)}, 
&\text{ if $2\,|\,m_{2i+1}(\mu),\ \forall i$}\\ 
0, &\text{ otherwise.}
\end{array} \right.
\end{split}
\end{equation*}
There are no similar formulas for $|\T_{(p^{+},p^{-})}(\mu)|$
or $|\T_{(p^{+},p^{-})}(\mu)^\psi|$.

\newpar
We label the irreducible characters of $S_{|\nu|}$ as
$\{\chi^{\rho}\,|\,\rho\vdash|\nu|\}$ as in \cite[I.7]{macdonald},
and write $\chi_\nu^{\rho}$ for the value of $\chi^{\rho}$ at an element
of cycle-type $\nu$,
so that $\chi^{(|\nu|)}$ is the trivial character and 
$\chi_\nu^{\rho'}=\epsilon_\nu\chi_\nu^{\rho}$.
\subsection{Character Theory of $GL_n(\Fq)$ and $U_n(\Fqs)$}
For the remainder of the paper, we specialize the context
of \S1.1 drastically, to the case
when $G\cong GL_n$, for some positive integer $n$. More concretely,
let $V$ be a vector space over $k$ of dimension $n$ and let $G=GL(V)$.
There are two kinds of $\Fq$-structures on $G$, split and non-split.
A Frobenius map $F:G\to G$ is split if it is induced by some
Frobenius map $F_V:V\to V$, in the sense that
\[ F_V(gv)=F(g)F_V(v),\ \forall g\in G, v\in V. \]
Then $G^{F}\cong \mathrm{Aut}_{\Fq}(V^{F_V})\cong GL_n(\Fq)$, the finite
general linear group. If $F$ is a non-split Frobenius map, there exists
some outer involution $\theta'$ of $G$ commuting with $F$, and for \emph{any}
such $\theta'$, $\theta' F$ is a split Frobenius map. In this case
$G^{F}\cong U_n(\Fqs)$, the finite unitary group.

\newpar
As in \cite[Chapter IV]{macdonald},
we will need to consider the
system of maps $\widehat{\Fqe^{\times}}\to
\widehat{\mathbb{F}_{q^{e'}}^{\times}}$ for $e\,|\,e'$
(the transpose of the norm map), and its limit
$L=\mathrm{colim}\ \widehat{\Fqe^{\times}}$.
Let $\sigma$ denote the $q$-th power map
on both $k^{\times}$ and $L$, so that 
$(k^{\times})^{\sigma^{e}}\cong\Fqe^{\times}$,
$L^{\sigma^{e}}\cong\widehat{\Fqe^{\times}}$ for all $e\geq 1$.
Write $\langle\cdot,\cdot\rangle^{\sigma^{e}}:
(k^{\times})^{\sigma^{e}}\times L^{\sigma^{e}} \to \Qlbar^{\times}$
for the canonical pairing.
Let $\iota$ denote the inverse map on $k^{\times}$ and $L$, and write
$\tsigma$ for $\iota\sigma$, the $(-q)$-th power map. (Note that
$\tsigma^{2}=\sigma^{2}$.) We also have a canonical pairing
$\langle\cdot,\cdot\rangle^{\tsigma^{e}}:
(k^{\times})^{\tsigma^{e}}\times L^{\tsigma^{e}} \to \Qlbar^{\times}$
(the same as $\langle\cdot,\cdot\rangle^{\sigma^{e}}$ if $e$ is even).

\newpar
We will fix some set of representatives for the orbits of the group
$\langle\sigma\rangle$ generated by $\sigma$ on $L$, and call it
$\langle\sigma\rangle\!\setminus\! L$. Similarly define
$\langle\tsigma\rangle\!\setminus\! L$. For
$\xi\in\langle\sigma\rangle\!\setminus\! L$, let
$m_\xi=|\langle\sigma\rangle.\xi|$, in other words the smallest $e\geq 1$
such that $\sigma^{e}(\xi)=\xi$. Let $d_\xi=
\langle -1,\xi \rangle^{\sigma^{m_\xi}}$, which equals
$1$ if $L^{\sigma^{m_\xi}}$ contains square roots of $\xi$, and $-1$ if it
does not.
Similarly define $\tilde{m}_\xi$ and $\tilde{d}_\xi$
using $\tsigma$ instead of $\sigma$.

\newpar
First consider the case when $F:G\to G$ is a split Frobenius map.
Let $\CalP_n$ be the set of collections of partitions
$\underline{\mu}=(\mu_\alpha)_{\alpha\in k^{\times}}$, almost all zero,
such that $\sum_{\alpha\in k^{\times}} |\mu_\alpha|= n$.
Let $\CalP_n^{\sigma}$ be the subset of $\CalP_n$ consisting of all
$\underline{\mu}$ such that
$\mu_{\sigma(\alpha)}=\mu_\alpha$ for all $\alpha$. It is well known
that there is a natural bijection between 
$\CalP_n^{\sigma}$ and the set of conjugacy classes in $G^{F}$.
Dually, let $\widehat{\CalP}_n$ be the set of collections
of partitions $\underline{\nu}=(\nu_\xi)_{\xi\in L}$, almost all zero,
such that $\sum_{\xi\in L} |\nu_\xi|= n$.
Let $\widehat{\CalP}_n^{\sigma}$ be the subset of $\widehat{\CalP}_n$
of all $\underline{\nu}$ such that $\nu_{\sigma(\xi)}=\nu_\xi$ for all $\xi$.
Note that for $\underline{\nu}\in\widehat{\CalP}_n^{\sigma}$,
\[ \sum_{\xi\in\langle\sigma\rangle\setminus L} m_\xi |\nu_\xi| = n. \]
For $\underline{\nu},\underline{\rho}\in\widehat{\CalP}_n^{\sigma}$,
we write $|\underline{\nu}|=|\underline{\rho}|$ to mean that
$|\nu_\xi|=|\rho_\xi|$ for all $\xi$.

\newpar
We can define a bijection
between $\widehat{\CalP}_n^{\sigma}$ and the set of $G^{F}$-orbits of pairs
$(T,\lambda)$ as in \S1.1, so that if $(T,\lambda)$ is in the orbit
corresponding to $\underline{\nu}$:
\begin{enumerate}
\item the eigenlines of $T$ can be labelled
\[ \{L_{(\xi,j,i)}\,|\, \xi\in\langle\sigma\rangle\!\setminus\! L,\
1\leq j \leq \ell(\nu_\xi),\ i\in\Z/m_\xi(\nu_\xi)_j\Z\} \]
so that under the resulting isomorphism
\[ T \cong \prod_{\xi\in\langle\sigma\rangle\setminus L}
\prod_{j=1}^{\ell(\nu_\xi)}
\underbrace{k^{\times}\times\cdots\times k^{\times}}
_{\text{$m_\xi(\nu_\xi)_j$ factors}}, \]
$F|_T$ corresponds to cyclic permutation
of each group of factors $k^{\times}$, composed with $\sigma$;
\item  consequently,
\[ T^{F} \cong \prod_{\xi\in\langle\sigma\rangle\setminus L}
\prod_{j=1}^{\ell(\nu_\xi)}
(k^{\times})^{\sigma^{m_\xi(\nu_\xi)_j}}; \]
\item  under this isomorphism, $\lambda$ corresponds to
\[ \prod_{\xi\in\langle\sigma\rangle\setminus L}
\prod_{j=1}^{\ell(\nu_\xi)} 
\langle\cdot,\xi\rangle^{\sigma^{m_\xi(\nu_\xi)_j}}. \]
\end{enumerate}

\newpar
For $\underline{\nu}\in\widehat{\CalP}_n^{\sigma}$, let
$B_{\underline{\nu}}=\mathrm{tr}(\cdot,R_T^{\lambda})$ for
$(T,\lambda)$ in the corresponding $G^{F}$-orbit. As proved by Lusztig
in \cite{coxeter},
these coincide with the \emph{basic characters} defined by Green 
in \cite{green}.
(In particular, their values are computable, but this is not relevant here.)
Green's main result on the character theory of $GL_n(\Fq)$ states that
for any $\underline{\rho}\in\widehat{\CalP}_n^{\sigma}$,
\[ \chi^{\underline{\rho}} := 
(-1)^{n+\sum_{\xi\in\langle\sigma\rangle\setminus L} |\rho_\xi|}
\sum_{\substack{\underline{\nu}\in\widehat{\CalP}_n^{\sigma}\\
|\underline{\nu}|=|\underline{\rho}|}}
\left(\prod_{\xi\in\langle\sigma\rangle\setminus L}
(z_{\nu_\xi})^{-1}\chi_{\nu_\xi}^{\rho_\xi}\right) B_{\underline{\nu}} \]
is an irreducible character of $G^{F}$, and all irreducible
characters arise in this way for unique
$\underline{\rho}\in\widehat{\CalP}_n^{\sigma}$.
(See also \cite[Chapter IV]{macdonald}
and \cite[Theorem 1.2.10]{bks}. Note that Macdonald's parameters in
$\widehat{\CalP}_n^{\sigma}$
differ from those of \cite{bks} by transposing all partitions; we are
following the convention of \cite{bks}.) 
In words, the transition matrix between
the basic characters and the irreducible characters is formed from
the character table of various symmetric groups. 
Inverting this matrix, we have that
for any $\underline{\nu}\in\widehat{\CalP}_n^{\sigma}$,
\begin{equation} \label{glneqn}
B_{\underline{\nu}}=
(-1)^{n+\sum_{\xi\in\langle\sigma\rangle\setminus L} |\nu_\xi|}
\sum_{\substack{\underline{\rho}\in\widehat{\CalP}_n^{\sigma}\\
|\underline{\rho}|=|\underline{\nu}|}}
\left(\prod_{\xi\in\langle\sigma\rangle\setminus L}
\chi_{\nu_\xi}^{\rho_\xi}\right) \chi^{\underline{\rho}}.
\end{equation}
There is an obvious action of $L^{\sigma}$ on $\widehat{\CalP}_n^{\sigma}$,
and in particular, for $\underline{\nu}\in\widehat{\CalP}_n^{\sigma}$,
$\zeta.\underline{\nu}$ is well defined. For any $\eta\in L^{\sigma}$,
$B_{\eta.\underline{\nu}}$ and $\chi^{\eta.\underline{\rho}}$ are the
result of multiplying $B_{\underline{\nu}}$ and $\chi^{\underline{\rho}}$
by the one-dimensional character $\langle\det(\cdot),\eta\rangle^{\sigma}$
of $G^{F}$. The \emph{unipotent} irreducible characters referred to in the
introduction are those $\chi^{\underline{\rho}}$ for which
$\rho_\xi=0$ unless $\xi=1$. (In the introduction we parametrized these
by $\rho=\rho_1$.)

\newpar
The case when $F:G\to G$ is a non-split Frobenius map is less well known,
but very similar, in fact mostly identical once $\sigma$ is replaced by
$\tsigma$, $m_\xi$ by $\tilde{m}_\xi$, and so on. Define
$\widehat{\CalP}_n^{\tsigma}$ in the obvious way. Again, for any
$\underline{\nu}\in\widehat{\CalP}_n^{\tsigma}$,
\[ \sum_{\xi\in\langle\tsigma\rangle\setminus L} \tilde{m}_\xi |\nu_\xi|= n. \]
For $\underline{\nu},\underline{\rho}\in\widehat{\CalP}_n^{\tsigma}$,
we write $|\underline{\nu}|=|\underline{\rho}|$ to mean that
$|\nu_\xi|=|\rho_\xi|$ for all $\xi$.

\newpar
We can define a bijection
between $\widehat{\CalP}_n^{\tsigma}$ and the set of $G^{F}$-orbits of pairs
$(T,\lambda)$ as above, so that if $(T,\lambda)$ is in the orbit
corresponding to $\underline{\nu}$:
\begin{enumerate}
\item the eigenlines of $T$ can be labelled
\[ \{L_{(\xi,j,i)}\,|\, \xi\in\langle\tsigma\rangle\!\setminus\! L,\
1\leq j \leq \ell(\nu_\xi),\ i\in\Z/\tilde{m}_\xi(\nu_\xi)_j\Z\} \]
so that under the resulting isomorphism
\[ T \cong \prod_{\xi\in\langle\tsigma\rangle\setminus L}
\prod_{j=1}^{\ell(\nu_\xi)}
\underbrace{k^{\times}\times\cdots\times k^{\times}}
_{\text{$\tilde{m}_\xi(\nu_\xi)_j$ factors}}, \]
$F|_T$ corresponds to cyclic permutation
of each group of factors $k^{\times}$, composed with $\tsigma$;
\item  consequently,
\[ T^{F} \cong \prod_{\xi\in\langle\tsigma\rangle\setminus L}
\prod_{j=1}^{\ell(\nu_\xi)}
(k^{\times})^{\tsigma^{\tilde{m}_\xi(\nu_\xi)_j}}; \]
\item  under this isomorphism, $\lambda$ corresponds to
\[ \prod_{\xi\in\langle\tsigma\rangle\setminus L}
\prod_{j=1}^{\ell(\nu_\xi)} 
\langle\cdot,\xi\rangle^{\tsigma^{\tilde{m}_\xi(\nu_\xi)_j}}. \]
\end{enumerate}

\newpar
For $\underline{\nu}\in\widehat{\CalP}_n^{\tsigma}$, let
$B_{\underline{\nu}}=\mathrm{tr}(\cdot,R_T^{\lambda})$ for
$(T,\lambda)$ in the corresponding $G^{F}$-orbit. The extension of Green's
result to the non-split case was proved by Lusztig and Srinivasan
in \cite[Theorem 3.2]{lusztigsrinivasan}: in our notation, for any
$\underline{\rho}\in\widehat{\CalP}_n^{\sigma}$,
\[ \chi^{\underline{\rho}} := 
(-1)^{\lceil\frac{n}{2}\rceil
+\sum_{\xi\in\langle\tsigma\rangle\setminus L} \tilde{m}_\xi n(\rho_\xi')
+ |\rho_\xi|}
\sum_{\substack{\underline{\nu}\in\widehat{\CalP}_n^{\tsigma}\\
|\underline{\nu}|=|\underline{\rho}|}}
\left(\prod_{\xi\in\langle\tsigma\rangle\setminus L}
(z_{\nu_\xi})^{-1}\chi_{\nu_\xi}^{\rho_\xi}\right) B_{\underline{\nu}} \]
is an irreducible character of $G^{F}$, and all irreducible
characters arise in this way for unique
$\underline{\rho}\in\widehat{\CalP}_n^{\tsigma}$. Inverting, we see that 
for any $\underline{\nu}\in\widehat{\CalP}_n^{\tsigma}$,
\begin{equation} \label{uneqn}
B_{\underline{\nu}}=
\sum_{\substack{\underline{\rho}\in\widehat{\CalP}_n^{\tsigma}\\
|\underline{\rho}|=|\underline{\nu}|}}
(-1)^{\lceil\frac{n}{2}\rceil
+\sum_{\xi\in\langle\tsigma\rangle\setminus L} \tilde{m}_\xi n(\rho_\xi')
+ |\rho_\xi|}
\left(\prod_{\xi\in\langle\tsigma\rangle\setminus L}
\chi_{\nu_\xi}^{\rho_\xi}\right) \chi^{\underline{\rho}}.
\end{equation}
Again, the obvious action of $L^{\tsigma}$ on $\widehat{\CalP}_n^{\tsigma}$
corresponds to multiplication by one-dimensional characters, and the
\emph{unipotent} irreducible characters of the introduction are
those $\chi^{\underline{\rho}}$ for which $\rho_\xi=0$ unless $\xi=1$.
(In contrast to the case of $GL_n(\Fq)$, not all the unipotent
characters are constituents of $\Ind_{B^{F}}^{G^{F}}(1)$ for an $F$-stable
Borel subgroup $B$.)
\subsection{Descriptions of Weyl Groups}
Much of this paper deals with the special properties of Weyl groups in
$GL_n$, so the following ideas and notation will be crucial.
Let $T$ be any maximal torus of $G=GL(V)$. If 
$\{L_\omega\,|\,\omega\in\Omega\}$ is some labelling of the eigenlines of $T$,
we can identify the Weyl group $W(T)$ with the 
group of permutations of $\Omega$. Let $W(T)_\inv$ be the set of involutions
in $W(T)$, and $W(T)_{\ff-\inv}$ the set of fixed-point free involutions
(note that this means fixed-point free as a permutation of $\Omega$,
not as an automorphism of $T$).

\newpar
If $T$ is $\theta$-stable, there is a special involution
$w_1^{T}\in W(T)_\inv$ characterized as follows.
If $\theta$ is an inner involution, namely conjugation by $s\in G$, then
\[ s(L_\omega)=L_{w_1^{T}(\omega)},\ \forall \omega\in\Omega, \text{ and }
\theta|_T=w_1^{T}.\]
If $\theta$ is an outer involution, namely adjoint inverse with respect
to a nondegenerate symplectic or symmetric form on $V$, then
\[ L_\omega^{\perp}=\bigoplus_{\omega'\neq w_1^{T}(\omega)}L_{\omega'},\ 
\forall \omega\in\Omega, \text{ and }
\theta|_T=w_1^{T}\circ\iota.\]
For any $T$, if $f\in\Theta_T$, then $f^{-1}Tf$ is a $\theta$-stable
maximal torus, and we obtain
\[ w_f^{T}=\mathrm{Ad}(f^{-1})\circ w_1^{f^{-1}Tf}\circ\mathrm{Ad}(f)\in
W(T)_\inv. \]
When $T$ is fixed, we will write $w_f$ instead of $w_f^{T}$.
Clearly $w_f$ depends only on the double coset $TfK$; the resulting
map $T\!\setminus\!\Theta_T/K\to W(T)_\inv$ will be crucial in our
combinatorial rewritings of Lusztig's theorem.

\newpar
Now suppose that $F:G\to G$ is split, and $(T,\lambda)$ is a pair in the
$G^{F}$-orbit corresponding to $\underline{\nu}\in\widehat{\CalP}_n^{\sigma}$.
Apart from the full Weyl group $W(T)$,
we will mainly be concerned with the $F$-fixed subgroup 
$W(T)^{F}$ and its subsets
\begin{equation*}
\begin{split}
W(T)^{F}_{\lambda}&=\{w\in W(T)^{F}\,|\,\lambda\circ w= \lambda\},\\
W(T)^{F}_{\lambda\to\lambda^{-1}}&=
\{w\in W(T)^{F}\,|\,\lambda\circ w= \lambda^{-1}\}.
\end{split}
\end{equation*}
We can consider these sets in the framework of \S1.2 as follows. Define 
\[ \Lambda(\underline{\nu})=\{(\xi,j)\in\langle\sigma\rangle\!\setminus\! L
\times\Z\,|\,1\leq j\leq \ell(\nu_\xi)\}. \]
By abuse of notation, write $\underline{\nu}$ also for the
partition of $n$ whose parts are 
\[ (m_\xi(\nu_\xi)_j)_{(\xi,j)\in\Lambda(\underline{\nu})}. \]
Then $\Lambda(\underline{\nu})$
indexes the parts of $\underline{\nu}$, in accordance with the notation
of \S1.2; also $\ell(\underline{\nu})=|\Lambda(\underline{\nu})|$
is the $\Fq$-rank of $T$.
Now $\Omega(\underline{\nu})$ is the set of triples
\[ \{(\xi,j,i)\,|\, (\xi,j)\in\Lambda(\underline{\nu}),
i\in\Z/m_\xi(\nu_\xi)_j\Z\}, \]
which is precisely the set of labels of the eigenlines of $T$ used in \S1.3.
So $W(T)$ in the above realization, namely as the group of permutations
of $\Omega(\underline{\nu})$, coincides with $S_{|\underline{\nu}|}$ in the 
realization of \S1.2. As in that section, let
$w_{\underline{\nu}}\in W(T)$
be the permutation $(\xi,j,i)\mapsto(\xi,j,i+1)$. 
Then by the description of $F|_T$ given in \S1.3,
$W(T)^{F}$ is exactly $Z_{W(T)}(w_{\underline{\nu}})=Z^{\underline{\nu}}$.
Note that
\[ \epsilon_{\underline{\nu}}=\mathrm{sign}(w_{\underline{\nu}})
=\prod_{\substack{\xi\in\langle\sigma\rangle\!\setminus\! L\\
2\nmid m_\xi}}\epsilon_{\nu_\xi}
\prod_{\substack{\xi\in\langle\sigma\rangle\!\setminus\! L\\
2|m_\xi}}(-1)^{\ell(\nu_\xi)}. \]

\newpar
For $w\in W(T)^{F}$, we will use the notation $\bar{w}$, 
$i(w,\xi,j)\in\Z/m_\xi(\nu_\xi)_j$ of \S1.2. So
\[ w(\xi,j,i)=(\bar{w}(\xi,j),i+i(w,\xi,j)),\ 
\forall (\xi,j)\in\Lambda(\underline{\nu}), i\in\Z/m_\xi(\nu_\xi)_j\Z. \]
With this notation, $w\in W(T)^{F}$ lies in $W(T)^{F}_{\lambda}$
if and only if $\bar{w}(\xi,j)=(\xi,j')$ (for
some $j'$) and $m_\xi\,|\,i(w,\xi,j)$ hold for all 
$(\xi,j)\in\Lambda(\underline{\nu})$.
Here it is helpful to consider 
$\coprod_{\xi\in\langle\sigma\rangle\!\setminus\! L} \Omega(\nu_\xi)$, which
is the set
of triples $(\xi,j,s)$ with $(\xi,j)\in\Lambda(\underline{\nu})$
and $s\in\Z/(\nu_\xi)_j\Z$, and
$\prod_{\xi\in\langle\sigma\rangle\!\setminus\! L} S_{|\nu_\xi|}$, which
is the group of permutations of such triples which preserve the first
factor. For $w\in W(T)^{F}_{\lambda}$, we define 
$\hat{w}\in\prod_{\xi\in\langle\sigma\rangle\!\setminus\! L} S_{|\nu_\xi|}$ by
\[ \hat{w}(\xi,j,s)=(\bar{w}(\xi,j),s+\frac{i(w,\xi,j)}{m_\xi}). \]
Clearly $w\mapsto \hat{w}$ is an isomorphism between 
$W(T)^{F}_{\lambda}$ and
$\prod_{\xi\in\langle\sigma\rangle\!\setminus\! L} Z^{\nu_\xi}$.

\newpar
The analogous description of $W(T)^{F}_{\lambda\to\lambda^{-1}}$ is as follows.
For $\xi\in\langle\sigma\rangle\!\setminus\! L$, let $\xi^{\vee}$
be the chosen representative in the $\langle\sigma\rangle$-orbit of $\xi^{-1}$.
Define $i_0(\xi)\in\Z/m_\xi\Z$ by 
\[ \xi^{\vee}=\sigma^{i_0(\xi)}(\xi^{-1}), \]
so that $i_0(\xi^{\vee})=-i_0(\xi)$.
Assume $\nu_{\xi^{\vee}}=\nu_\xi$ for all 
$\xi\in\langle\sigma\rangle\!\setminus\! L$; otherwise
$W(T)^{F}_{\lambda\to\lambda^{-1}}$ is empty.
Clearly $w\in W(T)^{F}$ lies in $W(T)^{F}_{\lambda\to\lambda^{-1}}$
iff $\bar{w}(\xi,j)=(\xi^{\vee},j')$ (for some $j'$)
and $i(w,\xi,j)\equiv i_0(\xi)\ \mathrm{mod}\ m_\xi$ hold for all 
$(\xi,j)\in\Lambda(\underline{\nu})$. For every
permutation $\bar{w}$ of $\Lambda(\underline{\nu})$
such that $\bar{w}(\xi,j)=(\xi^{\vee},j')$ (for some $j'$),
and every $(\xi,j)\in\Lambda(\underline{\nu})$, 
lift $i_0(\xi)$ to an element $i_0(\bar{w},\xi,j)$ of
$\Z/m_\xi(\nu_\xi)_j\Z$ in an arbitrary way.
Consider the set of permutations of
$\coprod_{\xi\in\langle\sigma\rangle\!\setminus\! L} \Omega(\nu_\xi)$
which interchange $\Omega(\nu_\xi)$ and $\Omega(\nu_{\xi^{\vee}})$ 
for all $\xi$.
Since $\nu_{\xi^{\vee}}=\nu_\xi$, we can identify this set with
$\prod_{\{\xi,\xi^{\vee}\}} S_{|\nu_\xi|}$.
For $w\in W(T)^{F}_{\lambda\to\lambda^{-1}}$, we define
$\tilde{w}\in\prod_{\{\xi,\xi^{\vee}\}} S_{|\nu_\xi|}$ by
\[ \tilde{w}(\xi,j,s)=(\bar{w}(\xi,j),s+
\frac{i(w,\xi,j)-i_0(\bar{w},\xi,j)}{m_\xi}). 
\]
Clearly $w\mapsto \tilde{w}$ defines a bijection between 
$W(T)^{F}_{\lambda\to\lambda^{-1}}$ and
$\prod_{\{\xi,\xi^{\vee}\}}Z^{\nu_\xi}$.

\newpar
Now consider involutions in $W(T)^{F}$. Since $W(T)^{F}_\inv
=Z^{\underline{\nu}}_\inv$, we can apply the concepts of \S1.2. 
In particular, any
$w\in W(T)^{F}_\inv$ decomposes $\Lambda(\underline{\nu})$ into
\begin{equation*}
\begin{split}
\Lambda_w^{1}(\underline{\nu})&=
\{(\xi,j)\in\Lambda(\underline{\nu})\,|\,
\bar{w}(\xi,j)=(\xi,j), i(w,\xi,j)=0\},\\
\Lambda_w^{2}(\underline{\nu})&=
\{(\xi,j)\in\Lambda(\underline{\nu})\,|\,
\bar{w}(\xi,j)=(\xi,j), i(w,\xi,j)=\frac{1}{2}m_\xi(\nu_\xi)_j\},\text{ and}\\
\Lambda_w^{3}(\underline{\nu})&=
\{(\xi,j)\in\Lambda(\underline{\nu})\,|\,
\bar{w}(\xi,j)\neq(\xi,j)\}.
\end{split}
\end{equation*}
If $w\in W(T)^{F}_{\ff-\inv}$, then $\Lambda_w^{1}(\underline{\nu})$ is empty.
As in \S1.2, we write $\ell_w^{i}(\underline{\nu})$ for
$|\Lambda_w^{i}(\underline{\nu})|$, and also use the notations
$\ell_w^{i}(\underline{\nu})_0$, $\ell_w^{i}(\underline{\nu})_1$,
$\ell_w^{i}(\underline{\nu})_{2\,\mathrm{mod}\,4}$ etc.
Note that
$W(T)^{F}_{\pm-\inv}$, i.e.\ the $F$-fixed points of
$W(T)_{\pm-\inv}=(S_{|\underline{\nu}|})_{\pm-\inv}$, is identified
with $Z^{\underline{\nu}}_{\pm-\inv}$ as defined in \S1.2, and similar
statements holds for $W(T)^{F}_{(p^{+},p^{-})-\inv}$ and
$W(T)^{F}_{\star-\inv}$.

\newpar
We will aso write $W(T)^{F}_{\lambda,\inv}$ for $W(T)^{F}_{\lambda}\cap 
W(T)^{F}_\inv$, and $W(T)^{F}_{\lambda\to\lambda^{-1},\inv}$
for $W(T)^{F}_{\lambda\to\lambda^{-1}}\cap W(T)^{F}_\inv$.
Under the isomorphism $W(T)^{F}_{\lambda}\isomto
\prod_{\xi} Z^{\nu_\xi}$ defined above, $W(T)^{F}_{\lambda,\inv}$
corresponds to $\prod_{\xi} Z^{\nu_\xi}_\inv$. Moreover, as will
be crucial later, if $w\in W(T)^{F}_{\lambda,\inv}$ corresponds to
$(w_\xi)$, then
\[ \ell_w^{1}(\underline{\nu})=
\sum_{\xi\in\langle\sigma\rangle\!\setminus\! L} \ell_{w_\xi}^{1}(\nu_\xi), \]
and similarly for $\ell^{2}$ and $\ell^{3}$. However, the situation
is more complex when we introduce divisibility criteria, because
of the division by $m_\xi$: for instance
\[ \ell_w^{1}(\underline{\nu})_0=
\sum_{\substack{\xi\in\langle\sigma\rangle\!\setminus\! L\\2\nmid m_\xi}}
\ell_{w_\xi}^{1}(\nu_\xi)_0
+\sum_{\substack{\xi\in\langle\sigma\rangle\!\setminus\! L\\2| m_\xi}}
\ell_{w_\xi}^{1}(\nu_\xi). \]

\newpar
To get a similar description of $W(T)^{F}_{\lambda\to\lambda^{-1},\inv}$,
we need to put further constraints on the choice of $i_0(\bar{w},\xi,j)$.
It is easy to see that we can arrange to have
\[ i_0(\bar{w},\bar{w}(\xi,j))=-i_0(\bar{w},\xi,j) \]
except in the case when $\xi^{\vee}=\xi$, $\xi\neq 1,-1$, $2|(\nu_\xi)_j$,
and $\bar{w}(\xi,j)=(\xi,j)$. (For instance, when
$\xi^{\vee}=\xi$, $\xi\neq 1,-1$, $2\nmid(\nu_\xi)_j$ and
$\bar{w}(\xi,j)=(\xi,j)$, we are forced to set
$i_0(\bar{w},\xi,j)=\frac{1}{2}m_\xi(\nu_\xi)_j$,
since $i_0(\xi)=\frac{1}{2}m_\xi$.)
But if $\bar{w}$ is induced as above from 
$w\in W(T)^{F}_{\lambda\to\lambda^{-1},\inv}$,
this case cannot arise, since if $2|(\nu_\xi)_j$ we cannot have both
\begin{equation*}
\begin{split}
i(w,\xi,j)&=-i(w,\xi,j)\ \mathrm{mod}\ m_\xi(\nu_\xi)_j, \text{ and}\\
i(w,\xi,j)&\equiv i_0(\xi)=\frac{1}{2}m_\xi\ \mathrm{mod}\ m_\xi. 
\end{split}
\end{equation*}
So $w\mapsto\tilde{w}$ with these conventions is a bijection between
$W(T)^{F}_{\lambda\to\lambda^{-1},\inv}$ and
\[ Z^{\nu_1}_\inv\times Z^{\nu_{-1}}_\inv\times
\prod_{\substack{\xi\in\langle\sigma\rangle\!\setminus\! L\\
\xi^{\vee}=\xi\\\xi\neq 1,-1}}
\{w_\xi\in Z^{\nu_\xi}_\inv\,|\,\ell_{w_\xi}^{1}(\nu_\xi)_0=
\ell_{w_\xi}^{2}(\nu_\xi)_0=0\}\times
\prod_{\substack{\{\xi_1\neq\xi_2\}\in\langle\sigma\rangle\!\setminus\! L\\
\xi_1^{\vee}=\xi_2}} Z^{\nu_{\xi_1}}_\inv. \]
In contrast to the previous situation, if we write $\tilde{w}$
as $(w_\xi)\in\prod_{\{\xi,\xi^{\vee}\}}Z^{\nu_\xi}_\inv$, then if
$\xi^{\vee}\neq\xi$, the elements of $\Lambda_{w_\xi}^{1}(\nu_\xi)$
contribute not to $\ell_w^{1}(\underline{\nu})$ but to
$\ell_w^{3}(\underline{\nu})$. Similarly, if
$\xi^{\vee}=\xi$, $\xi\neq 1,-1$, then $\Lambda_{w_\xi}^{1}(\nu_\xi)$
conributes to $\ell_w^{2}(\underline{\nu})$. This subtlety is one
factor complicating the formulas below.

\newpar
Everything we have said applies equally well to the non-split case, 
with $\tsigma$ instead of
$\sigma$ and $\tilde{m}_\xi$ instead of $m_\xi$ throughout, except
for the fact that the $\Fq$-rank of $T$ is not
$\ell(\underline{\nu})$ but $\ell(\underline{\nu})_0$.
Indeed, $q-1$ (as a polynomial) divides
\[ |(k^{\times})^{\tsigma^{\tilde{m}_\xi(\nu_\xi)_j}}|
=q^{\tilde{m}_\xi(\nu_\xi)_j}-(-1)^{\tilde{m}_\xi(\nu_\xi)_j} \]
once if $\tilde{m}_\xi(\nu_\xi)_j$ is even and not at all if
$\tilde{m}_\xi(\nu_\xi)_j$ is odd. In general, the main difference
between the split and non-split cases below lies in the calculation
of $\Fq$-ranks, which affects the signs in Lusztig's formula.
\section{Cases where $G/K=GL_n/Sp_n$}
In this section, we suppose that $V$ has a nondegenerate symplectic form
$\langle\cdot,\cdot\rangle$ (so in particular $n$ is even),
and that $\theta:G\to G$ is the involution defined by
\[ \langle\theta(g)v,v'\rangle=\langle v,g^{-1}v'\rangle,\ \forall g\in G,
v,v'\in V. \]
Since $G^{\theta}=Sp(V,\langle\cdot,\cdot\rangle)$ is connected, $K$
must be equal to it. So $G/K$ is the symmetric space $GL_n/Sp_n$.

\newpar
Let $T$ be any maximal torus of $G$. 
In \S1.4, we associated to any $f\in\Theta_T$ an involution
$w_f\in W(T)$, depending only on the double coset $TfK$.
Since $\langle\cdot,\cdot\rangle$ is a symplectic form, every line is
orthogonal to itself, so $w_f$ is fixed-point free (as a permutation
of the eigenlines of $T$, not as an automorphism of $T$). Thus
\[ T\cap fKf^{-1} = \{t\in T\,|\, w_f(t)=t^{-1}\} \]
is connected, and $Z_G(T\cap fKf^{-1})=T$. Moreover:
\begin{proposition} \label{glnspnprop}
The map $f\mapsto w_f$ induces a bijection $T\!\setminus\!\Theta_T/K
\isomto W(T)_{\ff-\inv}$.
\end{proposition}
\begin{proof}
This is very well known, especially when translated into the language
of flags via the connection mentioned in \S1.1. Surjectivity can be proved
by an explicit construction, and injectivity is easy by induction.
\end{proof}
\subsection{The $GL_n(\Fq)/Sp_n(\Fq)$ Case}
In this subsection, let $F:G\to G$ be a 
split Frobenius map which commutes with $\theta$.
So $F$ is induced by a Frobenius map $F_V$ on $V$ which respects
$\langle\cdot,\cdot\rangle$.
One has $G^{F}\cong GL_n(\Fq)$, $K^{F}\cong Sp_n(\Fq)$.
The following result was obtained by a different (and simpler)
method in \cite[\S4]{bks}:
\begin{theorem} \label{glnspnthm}
For any $\underline{\rho}\in\widehat{\CalP}_{n}^{\sigma}$,
\[ \langle \chi^{\underline{\rho}}, \Ind_{Sp_n(\Fq)}^{GL_n(\Fq)}(1)\rangle=
\left\{
\begin{array}{cl}
1, &\text{ if all $\rho_\xi$ are even}\\ 
0, &\text{ otherwise.} 
\end{array} \right. \]
\end{theorem}
By the results in \S1.3, it is equivalent to prove that for any
$\underline{\nu}\in\widehat{\CalP}_{n}^{\sigma}$,
\begin{equation} \label{glnspneqn}
\langle B_{\underline{\nu}}, \Ind_{K^{F}}^{G^{F}}(1)\rangle=
\prod_{\xi\in\langle\sigma\rangle\setminus L}
\sum_{\substack{\rho_\xi\vdash|\nu_\xi|\\
\rho_\xi\text{ even}}} \chi_{\nu_\xi}^{\rho_\xi}.
\end{equation}
Note that the sign in \eqref{glneqn} disappears since $n$ is even,
and $|\rho|$ is even for even $\rho$. 

\newpar
We aim to deduce this from
Lusztig's general formula for the left-hand side (Theorem \ref{lusztigthm}
above). This entails successively analysing the components of the
formula in our combinatorial terms, until we are reduced to a fact
about class functions on the symmetric group (in this case one which
is well known).
This basic strategy will be repeated in every case; the main reason for
including the present case, where the result is not new, is that it
serves as the prototype for the following subsections.

\newpar
Let $T$ be an $F$-stable maximal torus, and $\lambda:T^{F}\to\Qlbar^{\times}$
a character, for which $(T,\lambda)$ is in the $G^{F}$-orbit
corresponding to $\underline{\nu}\in\widehat{\CalP}_n^{\sigma}$.
We will use the description of $W(T)^{F}$ as $Z^{\underline{\nu}}$ given
in \S1.4. Proposition \ref{glnspnprop} implies:
\begin{lemma}
The map $f\mapsto w_f$ induces a bijection
\[ T^{F}\!\setminus\!\Theta_T^{F}/K^{F}\isomto
W(T)^{F}_{\ff-\inv}. \]
\end{lemma}
\begin{proof}
Clearly the map in Proposition \ref{glnspnprop} is $F$-stable. So we
need only note that since $T$, $K$, and all $T\cap fKf^{-1}$ are connected,
$T^{F}\!\setminus\!\Theta_T^{F}/K^{F}=(T\!\setminus\!\Theta_T/K)^{F}$.
\end{proof}
\begin{lemma} \label{glnspnlemma}
For $f\in\Theta_T^{F}$, $f\in\Theta_{T,\lambda}^{F}\Leftrightarrow
w_f\in W(T)_\lambda^{F}$.
\end{lemma}
\begin{proof}
Since $T\cap fKf^{-1}=\{t\in T\,|\,w_f(t)=t^{-1}\}$ is connected, 
$\epsilon_{T,f}=1$, so
\[ \Theta_{T,\lambda}^{F}=\{f\in\Theta_T^{F}\,|\, 
\lambda|_{\{t\in T^{F}|w_f(t)=t^{-1}\}}=1\}. \]
Thus it suffices to show that
\[ \{t\in T^{F}\,|\,w_f(t)=t^{-1}\}=\{tw_f(t)^{-1}\,|\,t\in T^{F}\}. \]
This follows easily from the fact that
$w_f\in W(T)^{F}_{\ff-\inv}$.
\end{proof}
\begin{corollary}
The map $f\mapsto w_f$ induces a bijection
\[ T^{F}\!\setminus\!\Theta_{T,\lambda}^{F}/K^{F}\isomto 
W(T)^{F}_{\lambda,\ff-\inv}. \]
\end{corollary}

\newpar
As noted above, $Z_G((T\cap fKf^{-1})^{\circ})=T$, so Lusztig's
formula becomes
\[ \langle B_{\underline{\nu}}, \Ind_{K^{F}}^{G^{F}}(1)\rangle=
|T^{F}\!\setminus\!\Theta_{T,\lambda}^{F}/K^{F}|=
|W(T)^{F}_{\lambda,\ff-\inv}|. \]
It is clear that under the isomorphism
$W(T)^{F}_{\lambda}\isomto
\prod_{\xi\in\langle\sigma\rangle\!\setminus\! L}Z^{\nu_\xi}$ given 
in \S1.4, $W(T)^{F}_{\lambda,\ff-\inv}$ corresponds to
$\prod_{\xi\in\langle\sigma\rangle\!\setminus\! L}Z^{\nu_\xi}_{\ff-\inv}$.
Now we apply the combinatorial fact 
(for which see \cite[VII.(2.4)]{macdonald}):
\begin{equation} \label{combglnspneqn}
|Z^{\nu}_{\ff-\inv}|
=\sum_{\substack{\rho\vdash|\nu|\\\rho\text{ even}}}
\chi_\nu^{\rho}.
\end{equation}
This gives Equation \eqref{glnspneqn} and hence Theorem \ref{glnspnthm}.
\subsection{The $U_n(\Fqs)/Sp_n(\Fq)$ Case}
Now we keep the assumptions from before \S2.1, but take $F$ to be instead a
non-split Frobenius map which commutes with $\theta$. So $\theta F$
is a split Frobenius map as above,
induced by a Frobenius map $F_V$ on $V$ which respects
$\langle\cdot,\cdot\rangle$. One has $G^{F}\cong U_n(\Fqs)$,
$K^{F}\cong Sp_n(\Fq)$. In this case the result is:
\begin{theorem} \label{unspnthm}
For any $\underline{\rho}\in\widehat{\CalP}_{n}^{\tsigma}$,
\[ \langle \chi^{\underline{\rho}}, \Ind_{Sp_n(\Fq)}^{U_n(\Fqs)}(1)\rangle=
\left\{
\begin{array}{cl}
1, &\text{ if all $\rho_\xi$ are even}\\ 
0, &\text{ otherwise.} 
\end{array} \right. \]
\end{theorem}
It is equivalent to prove that for any
$\underline{\nu}\in\widehat{\CalP}_{n}^{\tsigma}$,
\begin{equation} \label{unspneqn}
\langle B_{\underline{\nu}}, \Ind_{K^{F}}^{G^{F}}(1)\rangle=
\prod_{\xi\in\langle\tsigma\rangle\setminus L}
\sum_{\substack{\rho_\xi\vdash|\nu_\xi|\\
\rho_\xi\text{ even}}} \chi_{\nu_\xi}^{\rho_\xi}.
\end{equation}
Note that the sign in \eqref{uneqn} disappears because
$n(\rho')\equiv\frac{|\rho|}{2}$ mod $2$ for any
even $\rho$, so that
\[ \sum_{\xi\in\langle\tsigma\rangle\setminus L}
\tilde{m}_{\xi}n(\rho_\xi')\equiv
\sum_{\xi\in\langle\tsigma\rangle\setminus L}
\tilde{m}_{\xi}\frac{|\rho_\xi|}{2}\equiv\frac{n}{2} \text{ mod }2. \]
Equation \eqref{unspneqn} is deduced exactly as in the previous section,
using the non-split versions of the notation in \S1.4.
\section{Cases where $G/K=GL_n/(GL_{n^{+}}\times GL_{n^{-}})$}
\setcounter{theorem}{0}
In this section, we suppose that $\theta$ is an \emph{inner} involution, namely
conjugation with respect to $s\in G$ such that $s^{2}=1$.
Let $V^{+}$ be the $(+1)$-eigenspace and $V^{-}$ the $(-1)$-eigenspace
of $s$ on $V$, so that $V=V^{+}\oplus V^{-}$. Let $n^{+}=\dim V^{+}$,
$n^{-}=\dim V^{-}$, so that $n=n^{+}+n^{-}$. Since $G^{\theta}=
GL(V^{+})\times GL(V^{-})$ is connected, $K$ must be equal to it.
So $G/K$ is the symmetric space $GL_n/(GL_{n^{+}}\times GL_{n^{-}})$.

\newpar
Let $T$ be any maximal torus of $G$. For $f\in\Theta_T$,
we have $w_f\in W(T)_\inv$, not necessarily fixed-point free.
Indeed fixed points of $w_f$ correspond to eigenlines
of $T$ which are stable under $s$, and therefore lie in either $V^{+}$
or $V^{-}$. Let 
\[ \epsilon_f:\{\text{fixed points of $w_f$}\}\to\{+,-\} \]
be the resulting map. 
In the notation of \S1.2, $(w_f,\epsilon_f)\in W(T)_{(n^{+},n^{-})-\inv}$.
\begin{proposition} \label{glnglnglnprop}
The map $f\mapsto (w_f,\epsilon_f)$ induces a bijection
\[ T\!\setminus\!\Theta_T/K\isomto W(T)_{(n^{+},n^{-})-\inv}. \]
\end{proposition}
\begin{proof}
As with Proposition \ref{glnspnprop}, this is well known 
when translated in terms of $K$-orbits on the flag variety,
and easy to prove (see for instance \cite{yamamoto}).
\end{proof}
\subsection{The $GL_n(\Fq)/(GL_{n^{+}}(\Fq)\times GL_{n^{-}}(\Fq))$ Case}
In this subsection, assume that $F:G\to G$ is a split Frobenius map
such that $F(s)=s$. So $F$ is induced by a Frobenius map $F_V$ on $V$
which stabilizes $V^{+}$ and $V^{-}$. One has $G^{F}\cong GL_n(\Fq)$,
$K^{F}\cong GL_{n^{+}}(\Fq)\times GL_{n^{-}}(\Fq)$. 
Recall the definition of $\T_{(p^{+},p^{-})}(\mu)$ from \S1.2.
\begin{theorem} \label{glnglnglnthm}
For any $\underline{\rho}\in\widehat{\CalP}_{n}^{\sigma}$,
\[ \langle \chi^{\underline{\rho}},
\Ind_{GL_{n^{+}}(\Fq)\times GL_{n^{-}}(\Fq)}^{GL_n(\Fq)}(1)\rangle=
\left\{
\begin{array}{cl}
|\T_{(n^{+},n^{-})}(\rho_1')|, &\begin{array}{l}
\text{if $\rho_\xi=\rho_{\xi^{-1}}, \forall\xi$}\\
\text{and $\rho_{-1}'$ is even}\end{array}\\ 
0, &\text{ otherwise.} 
\end{array} \right. \]
\end{theorem}
\noindent
(Note that $|\T_{(n^{+},n^{-})}(\rho_1')|$ could be zero.)
By \eqref{glneqn}, it is equivalent to say that for any
$\underline{\nu}\in\widehat{\CalP}_{n}^{\sigma}$,
\begin{equation} \label{glnglnglneqn}
\begin{split}
\langle B_{\underline{\nu}}, \Ind_{K^{F}}^{G^{F}}(1)\rangle&=
(\sum_{\rho_1\vdash|\nu_1|}|\T_{(n^{+},n^{-})}(\rho_1')|\chi_{\nu_1}^{\rho_1})
(\sum_{\substack{\rho_{-1}\vdash|\nu_{-1}|\\\text{$\rho_{-1}'$ even}}}
\chi_{\nu_{-1}}^{\rho_{-1}})\\
&\qquad\times
\prod_{\substack{\xi\in\langle\sigma\rangle\setminus L\\
\xi^{\vee}=\xi\\ \xi\neq 1, -1}}
((-1)^{|\nu_\xi|}\sum_{\rho_\xi\vdash|\nu_\xi|} \chi_{\nu_\xi}^{\rho_\xi})
\prod_{\substack{\{\xi_1\neq\xi_2\}\in\langle\sigma\rangle\setminus L\\
\xi_1^{\vee}=\xi_2}} 
\delta_{\nu_{\xi_1}\nu_{\xi_2}}z_{\nu_{\xi_1}}.
\end{split}
\end{equation}
Here the sign $(-1)^{n+\sum_{\xi\in\langle\sigma\rangle\!\setminus\! L}
|\nu_\xi|}=(-1)^{\sum_{\xi\in\langle\sigma\rangle\!\setminus\! L}
(m_\xi+1)|\nu_\xi|}$ in \eqref{glneqn} is 
simplified by noting that $m_1=m_{-1}=1$,
$m_\xi$ is even for $\xi\neq 1, -1$ such that 
$\xi^{\vee}=\xi$, and $\nu_{\xi_1}$ and $\nu_{\xi_2}$
must be equal whenever $\xi_1^{\vee}=\xi_2$, 
$\xi_1\neq \xi_2$ in order for the right-hand side to be nonzero.

\newpar
Let $(T,\lambda)$ be a pair in the $G^{F}$-orbit corresponding to
$\underline{\nu}\in\widehat{\CalP}_{n}^{\sigma}$.
\begin{lemma} \label{trivglnglnglnlemma}
The map $f\mapsto (w_f,\epsilon_f)$ induces a bijection
\[ T^{F}\!\setminus\!\Theta_T^{F}/K^{F}\isomto
W(T)^{F}_{(n^{+},n^{-})-\inv}. \]
\end{lemma}
\begin{proof}
Again, this follows from Proposition \ref{glnglnglnprop}, since
$T$, $K$, and $T\cap fKf^{-1}=T^{w_f}$ are all connected.
\end{proof}
Note that for $w\in W(T)^{F}_{\lambda\to\lambda^{-1},\inv}$,
$(\xi,j)\in\Lambda_w^{1}(\underline{\nu})\Rightarrow \xi=\pm 1$. Define
\[ Y^{\underline{\nu}}_\inv=
\{w\in W(T)^{F}_{\lambda\to\lambda^{-1},\inv}\,|\,
(\xi,j)\in\Lambda_w^{1}(\underline{\nu})\Rightarrow \xi=1\}. \]
\begin{lemma} \label{glnglnglnlemma}
For $f\in\Theta_T^{F}$, $f\in\Theta_{T,\lambda}^{F}\Leftrightarrow
w_f\in Y^{\underline{\nu}}_\inv$.
\end{lemma}
\begin{proof}
Since $T\cap fKf^{-1}=T^{w_f}$ is connected, $\epsilon_{T,f}=1$, so
\[ \Theta_{T,\lambda}^{F}=\{f\in\Theta_T^{F}\,|\, 
\lambda|_{(T^{F})^{w_f}}=1\}. \]
Thus it suffices to show that $(T^{F})^{w_f}$ is generated by
\[ \{tw_f(t)\,|\,t\in T^{F}\} \text{ and }
\prod_{(\xi,j)\in\Lambda_{w_f}^{1}(\underline{\nu})}
(k^{\times})^{\sigma^{m_\xi(\nu_\xi)_j}}. \]
This follows easily from our
description of $T^{F}$.
\end{proof}
\begin{corollary}
The map $f\mapsto (w_f,\epsilon_f)$ induces a bijection
\[ T^{F}\!\setminus\!\Theta_{T,\lambda}^{F}/K^{F}\isomto 
Y^{\underline{\nu}}_{(n^{+},n^{-})-\inv}, \]
where $Y^{\underline{\nu}}_{(n^{+},n^{-})-\inv}=\{(w,\epsilon)
\in Z^{\underline{\nu}}_{(n^{+},n^{-})-\inv}\,|\,w\in Y^{\underline{\nu}}_\inv\}$.
\end{corollary}

\newpar
Now the $\Fq$-rank of $T$ is $\ell(\underline{\nu})$,
and that of
\begin{equation*}
\begin{split}
Z_G((T\cap fKf^{-1})^{\circ})
&=\prod_{\substack{(\xi,j,i)\\(\xi,j)\in\Lambda_{w_f}^{1}(\underline{\nu})}}
GL(L_{(\xi,j,i)})\\
&\thickspace\times
\prod_{\substack{\{(\xi,j,i),(\xi,j,i+\frac{1}{2}m_\xi(\nu_\xi)_j)\}\\
(\xi,j)\in\Lambda_{w_f}^{2}(\underline{\nu})}}
GL(L_{(\xi,j,i)}\oplus L_{(\xi,j,i+\frac{1}{2}m_\xi(\nu_\xi)_j)})\\
&\thickspace\times\prod_{\substack{\{(\xi,j,i),w_f(\xi,j,i)\}\\
(\xi,j)\in\Lambda_{w_f}^{3}(\underline{\nu})}}
GL(L_{(\xi,j,i)}\oplus L_{w_f(\xi,j,i)}) 
\end{split}
\end{equation*}
is $\ell(\underline{\nu})+\ell_{w_f}^{2}(\underline{\nu})$.
So Lusztig's formula becomes
\[
\langle B_{\underline{\nu}}, \Ind_{K^{F}}^{G^{F}}(1)\rangle
=\sum_{f\in T^{F}\setminus\Theta_{T,\lambda}^{F}/K^{F}}
(-1)^{\ell_{w_f}^{2}(\underline{\nu})}
=\negthickspace\negthickspace
\sum_{(w,\epsilon)\in Y^{\underline{\nu}}_{(n^{+},n^{-})-\inv}}
\negthickspace\negthickspace\negthickspace
(-1)^{\ell_{w}^{2}(\underline{\nu})}. \]

\newpar
Under the bijection $w\mapsto\tilde{w}$ defined in \S1.4,
$Y^{\underline{\nu}}_\inv$ corresponds to
\[ Z^{\nu_1}_\inv\times Z^{\nu_{-1}}_{\ff-\inv}\times
\prod_{\substack{\xi\in\langle\sigma\rangle\!\setminus\! L\\
\xi^{\vee}=\xi\\\xi\neq 1,-1}}
\{w_\xi\in Z^{\nu_\xi}_\inv\,|\,\ell_{w_\xi}^{1}(\nu_\xi)_0=
\ell_{w_\xi}^{2}(\nu_\xi)_0=0\}\times
\prod_{\substack{\{\xi_1\neq\xi_2\}\in\langle\sigma\rangle\!\setminus\! L\\
\xi_1^{\vee}=\xi_2}} Z^{\nu_{\xi_1}}_\inv. \]
Hence
\begin{equation*}
\begin{split}
\langle B_{\underline{\nu}}, \Ind_{K^{F}}^{G^{F}}(1)\rangle
&=\left(\sum_{(w_1,\epsilon_1)\in Z^{\nu_1}_{(n^{+},n^{-})-\inv}}
\negthickspace\negthickspace\negthickspace
(-1)^{\ell_{w_1}^{2}(\nu_1)}\right)
\left(\epsilon_{\nu_{-1}}|Z^{\nu_{-1}}_{\ff-\inv}|\right)\\
&\thickspace\times\prod_{\substack{\xi\in\langle\sigma\rangle\setminus L\\
\xi^{\vee}=\xi\\ \xi\neq 1, -1}}
(-1)^{|\nu_\xi|}|\{w_\xi\in Z^{\nu_\xi}_{\inv}\,|\,
\ell_{w_\xi}^{1}(\nu_\xi)_0=\ell_{w_\xi}^{2}(\nu_\xi)_0=0\}|\\
&\thickspace\times
\prod_{\substack{\{\xi_1\neq\xi_2\}\in\langle\sigma\rangle\setminus L\\
\xi_1^{\vee}=\xi_2}} 
\delta_{\nu_{\xi_1}\nu_{\xi_2}}z_{\nu_{\xi_1}}.
\end{split}
\end{equation*}
In the second factor we have used the fact that if 
$w_{-1}\in Z^{\nu_{-1}}_{\ff-\inv}$, then
\[ \ell_{w_{-1}}^{2}(\nu_{-1})
= \ell_{w_{-1}}^{2}(\nu_{-1})_0
\equiv \ell(\nu_{-1})_0 \ \mathrm{mod}\ 2. \]
In the third factor, we have noted that both $\Lambda_{w_\xi}^{1}(\nu_\xi)$ and
$\Lambda_{w_\xi}^{2}(\nu_\xi)$ contribute to 
$\ell_w^{2}(\underline{\nu})$, and if there exists
$w_\xi\in Z^{\nu_\xi}$ such that
$\ell_{w_\xi}^{1}(\nu_\xi)_0=\ell_{w_\xi}^{2}(\nu_\xi)_0=0$, then
\[ \ell_{w_\xi}^{1}(\nu_\xi)+\ell_{w_\xi}^{2}(\nu_\xi)
\equiv \ell(\nu_\xi)_1 \equiv |\nu_\xi|\ \mathrm{mod}\ 2.\]
In the fourth factor there is no contribution to 
$\ell_w^{2}(\underline{\nu})$, as noted in \S1.4.
Now we apply \eqref{combglnspneqn} and 
the following combinatorial facts (for which see \S5):
\begin{gather}
\label{combglnglnglneqn}
\sum_{(w,\epsilon)\in Z^{\nu}_{(p^{+},p^{-})-\inv}}\negthickspace\negthickspace
(-1)^{\ell_w^{2}(\nu)}
=\sum_{\rho\vdash|\nu|} |\T_{(p^{+},p^{-})}(\rho')|\chi_\nu^{\rho},\\
\label{othercombglnoneqn}
|\{w\in Z^{\nu}_{\inv}\,|\,
\ell_{w}^{1}(\nu)_0=\ell_{w}^{2}(\nu)_0=0\}|
=\sum_{\rho\vdash|\nu|}\chi_\nu^{\rho}. 
\end{gather}
These give Equation \eqref{glnglnglneqn} and hence Theorem \ref{glnglnglnthm}.
\subsection{The $GL_n(\Fq)/GL_{n/2}(\Fqs)$ Case}
Now keep the assumptions from before \S3.1, 
but take $F$ to be a split Frobenius map
such that $F(s)=-s$. (Then $F$ still commutes with $\theta$.)
So $F$ is induced by a Frobenius map $F_V$ on $V$ which interchanges
$V^{+}$ and $V^{-}$, whence $n$ is even and $n^{+}=n^{-}=\frac{n}{2}$.
One has $G^{F}\cong GL_n(\Fq)$, $K^{F}\cong GL_{n/2}(\Fqs)$. The result is:
\begin{theorem} \label{glnglnthm}
For any $\underline{\rho}\in\widehat{\CalP}_{n}^{\sigma}$,
\[ \langle \chi^{\underline{\rho}},
\Ind_{GL_{n/2}(\Fqs)}^{GL_n(\Fq)}(1)\rangle=
\left\{
\begin{array}{cl}
1, &
\text{if $\rho_\xi=\rho_{\xi^{-1}}, \forall\xi$,
$\rho_{1}$ is even, and $\rho_{-1}'$ is even}\\ 
0, &\text{otherwise.} 
\end{array} \right. \]
\end{theorem}
By \eqref{glneqn}, it is equivalent to say that for any
$\underline{\nu}\in\widehat{\CalP}_{n}^{\sigma}$,
\begin{equation} \label{glnglneqn}
\begin{split}
\langle B_{\underline{\nu}}, \Ind_{K^{F}}^{G^{F}}(1)\rangle&=
(\sum_{\substack{\rho_1\vdash|\nu_1|\\\text{$\rho_{1}$ even}}}
\chi_{\nu_1}^{\rho_1})
(\sum_{\substack{\rho_{-1}\vdash|\nu_{-1}|\\\text{$\rho_{-1}'$ even}}}
\chi_{\nu_{-1}}^{\rho_{-1}})\\
&\qquad\times
\prod_{\substack{\xi\in\langle\sigma\rangle\setminus L\\
\xi^{\vee}=\xi\\ \xi\neq 1, -1}}
((-1)^{|\nu_\xi|}\sum_{\rho_\xi\vdash|\nu_\xi|} \chi_{\nu_\xi}^{\rho_\xi})
\prod_{\substack{\{\xi_1\neq\xi_2\}\in\langle\sigma\rangle\setminus L\\
\xi_1^{\vee}=\xi_2}} 
\delta_{\nu_{\xi_1}\nu_{\xi_2}}z_{\nu_{\xi_1}}.
\end{split}
\end{equation}
(For the signs here, see the comments after Equation \eqref{glnglnglneqn}.)

\newpar
Let $(T,\lambda)$ be a pair in the $G^{F}$-orbit corresponding to
$\underline{\nu}\in\widehat{\CalP}_{n}^{\sigma}$. Obviously 
Lemma \ref{trivglnglnglnlemma} becomes:
\begin{lemma}
The map $f\mapsto (w_f,\epsilon_f)$ induces a bijection
\[ T^{F}\!\setminus\!\Theta_T^{F}/K^{F}\isomto
W(T)^{F}_{\star-\inv}. \]
\end{lemma}
Now Lemma \ref{glnglnglnlemma} holds again here, with the same 
$Y^{\underline{\nu}}_\inv$
and the same proof. So arguing as in \S3.1, we get
\[ \langle B_{\underline{\nu}}, \Ind_{K^{F}}^{G^{F}}(1)\rangle
=\negthickspace\negthickspace
\sum_{(w,\epsilon)\in Y^{\underline{\nu}}_{\star-\inv}}
\negthickspace\negthickspace\negthickspace\negthickspace
(-1)^{\ell_w^{2}(\underline{\nu})}, \]
where the definition of $Y^{\underline{\nu}}_{\star-\inv}$ is the obvious one.
The rest of the proof is also the same as in \S3.1, except that Equation
\eqref{combglnglnglneqn} is replaced by:
\begin{equation} \label{combglnglneqn}\negthickspace\negthickspace
\sum_{(w,\epsilon)\in Z^{\nu}_{\star-\inv}}\negthickspace\negthickspace
\negthickspace\negthickspace
(-1)^{\ell_w^{2}(\nu)}
=\sum_{\substack{\rho\vdash|\nu|\\\text{$\rho$ even}}}
\chi_\nu^{\rho}.
\end{equation}
This too will be proved in \S5.
\subsection{The $U_n(\Fqs)/(U_{n^{+}}(\Fqs)\times U_{n^{-}}(\Fqs))$ Case}
Still under the general assumptions of this section, 
let $F:G\to G$ be a non-split
Frobenius map for which $F(s)=s$. Replacing $s$ by a $G^{F}$-conjugate
if necessary, we may assume that there is some nondegenerate
symmetric form $\langle\cdot,\cdot\rangle$ on $V$, for which
$V^{+}$ and $V^{-}$ are orthogonal, and such that the associated
outer involution $\theta':G\to G$ commutes with $F$. Then
$\theta'F$ is the split Frobenius map induced by some $F_V:V\to V$
which respects $\langle\cdot,\cdot\rangle$ and fixes $V^{+}$ and $V^{-}$.
One has $G^{F}\cong U_n(\Fqs)$, 
$K^{F}\cong U_{n^{+}}(\Fqs)\times U_{n^{-}}(\Fqs)$.
\begin{theorem} \label{unununthm}
For any $\underline{\rho}\in\widehat{\CalP}_{n}^{\tsigma}$,
\[ \langle \chi^{\underline{\rho}},
\Ind_{U_{n^{+}}(\Fqs)\times U_{n^{-}}(\Fqs)}^{U_n(\Fqs)}(1)\rangle=
\left\{
\begin{array}{cl}
|\T_{(n^{+},n^{-})}(\rho_1')^{\psi}|, &\begin{array}{l}
\text{if $\rho_\xi=\rho_{\xi^{-1}}, \forall\xi$}\\
\text{and $\rho_{-1}'$ is even}\end{array}\\ 
0, &\text{ otherwise.} 
\end{array} \right. \]
\end{theorem}
\noindent
(Note that $|\T_{(n^{+},n^{-})}(\rho_1')^{\psi}|$ could be zero.)
By \eqref{uneqn}, it is equivalent to say that for any
$\underline{\nu}\in\widehat{\CalP}_{n}^{\tsigma}$,
\begin{equation} \label{unununeqn}
\begin{split}
\langle B_{\underline{\nu}}, \Ind_{K^{F}}^{G^{F}}(1)\rangle&=
(\sum_{\rho_1\vdash|\nu_1|}(-1)^{n(\rho_1)}
|\T_{(n^{+},n^{-})}(\rho_1')^{\psi}|\chi_{\nu_1}^{\rho_1})
((-1)^{\frac{|\nu_{-1}|}{2}}
\sum_{\substack{\rho_{-1}\vdash|\nu_{-1}|\\\text{$\rho_{-1}'$ even}}}
\chi_{\nu_{-1}}^{\rho_{-1}})\\
&\qquad\times\prod_{\substack{\xi\in\langle\tsigma\rangle\setminus L\\
\xi^{\vee}=\xi\\ \xi\neq 1, -1\\ 4|\tilde{m}_\xi}}
((-1)^{|\nu_\xi|}\sum_{\rho_\xi\vdash|\nu_\xi|} \chi_{\nu_\xi}^{\rho_\xi})
\prod_{\substack{\xi\in\langle\tsigma\rangle\setminus L\\
\xi^{\vee}=\xi\\ \xi\neq 1, -1\\ 4\nmid\tilde{m}_\xi}}
(\sum_{\rho_\xi\vdash|\nu_\xi|} \chi_{\nu_\xi}^{\rho_\xi})\\
&\qquad\times
\prod_{\substack{\{\xi_1\neq\xi_2\}\in\langle\tsigma\rangle\setminus L\\
\xi_1^{\vee}=\xi_2\\ 2|\tilde{m}_{\xi_1}}} 
\delta_{\nu_{\xi_1}\nu_{\xi_2}}z_{\nu_{\xi_1}}
\prod_{\substack{\{\xi_1\neq\xi_2\}\in\langle\tsigma\rangle\setminus L\\
\xi_1^{\vee}=\xi_2\\ 2\nmid\tilde{m}_{\xi_1}}} 
(-1)^{|\nu_{\xi_1}|}\delta_{\nu_{\xi_1}\nu_{\xi_2}}z_{\nu_{\xi_1}}.
\end{split}
\end{equation}

Here the sign $(-1)^{\lceil\frac{n}{2}\rceil+
\sum_{\xi\in\langle\tsigma\rangle\setminus L} \tilde{m}_\xi n(\rho_\xi')
+|\rho_\xi|}$ of \eqref{uneqn} is simplified as follows.
If $\rho_\xi=\rho_{\xi^{-1}}$ for all $\xi$
and $\rho_{-1}'$ is even, then
\[ n=\sum_{\xi\in\langle\tsigma\rangle\setminus L}\tilde{m}_\xi |\rho_\xi|
\equiv |\rho_1| \text{ mod } 2, \]
and the sign can be replaced by
\begin{equation*}
\begin{split}
(-1)^{\lfloor\frac{|\rho_1|}{2}\rfloor+n(\rho_1')}
&(-1)^{\frac{|\rho_{-1}|}{2}+n(\rho_{-1}')}
\prod_{\substack{\xi\in\langle\tsigma\rangle\setminus L\\
\xi^{\vee}=\xi\\ \xi\neq 1, -1}}
(-1)^{\frac{\tilde{m}_\xi}{2}(|\rho_\xi|+2n(\rho_\xi'))+|\rho_\xi|}\\
&\qquad\times
\prod_{\substack{\{\xi_1\neq\xi_2\}\in\langle\tsigma\rangle\setminus L\\
\xi_1^{\vee}=\xi_2}}
(-1)^{\tilde{m}_{\xi_1}(|\rho_{\xi_1}|+2n(\rho_{\xi_1}'))+2|\rho_{\xi_1}|}.
\end{split}
\end{equation*}
Then we observe that if $\rho_{-1}'$ is even, $n(\rho_{-1}')$ is even;
and if all even parts of $\rho_1'$ occur with even multiplicity
(which is necessary for $\bar{a}(n^{+},n^{-},\rho_1')\neq 0$), then
the Young diagram of $\rho_1'$ (excluding the top left corner if
$|\rho_1|$ is odd) can be tiled by $2\times 1$ dominoes, from which we
see that 
$\lfloor\frac{|\rho_1|}{2}\rfloor+n(\rho_1')+n(\rho_1)$ is even.

\newpar
Let $(T,\lambda)$ be a pair in the $G^{F}$-orbit corresponding to
$\underline{\nu}\in\widehat{\CalP}_{n}^{\tsigma}$. Proposition
\ref{glnglnglnprop} has the following corollary in this case:
\begin{lemma}
The map $f\mapsto (w_f,\epsilon_f)$ induces a bijection
\[ T^{F}\!\setminus\!\Theta_T^{F}/K^{F}\isomto
W(T)^{F}_{(n^{+},n^{-})-\inv}. \]
\end{lemma}
Define $Y^{\underline{\nu}}_\inv\subset 
W(T)^{F}_{\lambda\to\lambda^{-1},\inv}$ 
in the same way as in \S3.1, but in the $\tsigma$ version.
\begin{lemma} \label{unununlemma}
For $f\in\Theta_T^{F}$, $f\in\Theta_{T,\lambda}^{F}\Leftrightarrow
w_f\in Y^{\underline{\nu}}_\inv$.
\end{lemma}
\begin{proof}
The proof is exactly analogous to that of Lemma \ref{glnglnglnlemma}.
\end{proof}
\begin{corollary}
The map $f\mapsto (w_f,\epsilon_f)$ induces a bijection
\[ T^{F}\!\setminus\!\Theta_{T,\lambda}^{F}/K^{F}\isomto 
Y^{\underline{\nu}}_{(n^{+},n^{-})-\inv}. \]
\end{corollary}

\newpar
The important point of difference from \S3.1 is the $\Fq$-ranks involved.
The $\Fq$-rank of $T$ is now $\ell(\underline{\nu})_0$,
and that of
\begin{equation*}
\begin{split}
Z_G((T\cap fKf^{-1})^{\circ})
&=\prod_{\substack{(\xi,j,i)\\(\xi,j)\in\Lambda_{w_f}^{1}(\underline{\nu})}}
GL(L_{(\xi,j,i)})\\
&\thickspace\times
\prod_{\substack{\{(\xi,j,i),(\xi,j,i+\frac{1}{2}\tilde{m}_\xi(\nu_\xi)_j)\}\\
(\xi,j)\in\Lambda_{w_f}^{2}(\underline{\nu})}}
GL(L_{(\xi,j,i)}\oplus L_{(\xi,j,i+\frac{1}{2}\tilde{m}_\xi(\nu_\xi)_j)})\\
&\thickspace\times\prod_{\substack{\{(\xi,j,i),w_f(\xi,j,i)\}\\
(\xi,j)\in\Lambda_{w_f}^{3}(\underline{\nu})}}
GL(L_{(\xi,j,i)}\oplus L_{w_f(\xi,j,i)}) 
\end{split}
\end{equation*}
is $\ell(\underline{\nu})_0 + 
\ell_{w_f}^{2}(\underline{\nu})_{0\,\mathrm{mod}\,4} +
\frac{1}{2}\ell_{w_f}^{3}(\underline{\nu})_1$.
(The $GL_2$ factors corresponding to $(\xi,j)\in\Lambda_{w_f}^{2}
(\underline{\nu})$ are split iff
$4|\tilde{m}_\xi(\nu_\xi)_j$, and those corresponding to
$\{(\xi,j)\neq\bar{w}_f(\xi,j)\}$ are split iff
$2|\tilde{m}_\xi(\nu_\xi)_j$.)
So Lusztig's formula becomes
\begin{equation*}
\begin{split}
\langle B_{\underline{\nu}}, \Ind_{K^{F}}^{G^{F}}(1)\rangle
&=\sum_{f\in T^{F}\setminus\Theta_{T,\lambda}^{F}/K^{F}}
(-1)^{\ell_{w_f}^{2}(\underline{\nu})_{0\,\mathrm{mod}\,4} +
\frac{1}{2}\ell_{w_f}^{3}(\underline{\nu})_1}\\
&=\negthickspace\negthickspace\negthickspace\negthickspace\negthickspace
\sum_{(w,\epsilon)\in Y^{\underline{\nu}}_{(n^{+},n^{-})-\inv}}
\negthickspace\negthickspace\negthickspace\negthickspace\negthickspace
(-1)^{\ell_{w}^{2}(\underline{\nu})_{0\,\mathrm{mod}\,4} +
\frac{1}{2}\ell_{w}^{3}(\underline{\nu})_1}.
\end{split}
\end{equation*}

\newpar
Using the same reasoning as in \S3.1, we can transform this expression to get:
\begin{equation*}
\begin{split}
\langle B_{\underline{\nu}}, \Ind_{K^{F}}^{G^{F}}(1)\rangle
&=\negthickspace\negthickspace
\sum_{(w_1,\epsilon_1)\in Z^{\nu_1}_{(n^{+},n^{-})-\inv}}
\negthickspace\negthickspace
(-1)^{\ell_{w_1}^{2}(\nu_1)_{0\,\mathrm{mod}\,4}
+\frac{1}{2}\ell_{w_1}^{3}(\nu_1)_1}
\\
&\thickspace\times(-1)^{\frac{|\nu_{-1}|}{2}}\epsilon_{\nu_{-1}}
|Z^{\nu_{-1}}_{\ff-\inv}|\\
&\thickspace\times\prod_{\substack{\xi\in\langle\tsigma\rangle\setminus L\\
\xi^{\vee}=\xi\\ \xi\neq 1, -1\\4|\tilde{m}_\xi}}
(-1)^{|\nu_\xi|}|\{w_\xi\in Z^{\nu_\xi}_\inv\,|\,
\ell_{w_\xi}^{1}(\nu_\xi)=\ell_{w_\xi}^{2}(\nu_\xi)=0\}|\\
&\thickspace\times\prod_{\substack{\xi\in\langle\tsigma\rangle\setminus L\\
\xi^{\vee}=\xi\\ \xi\neq 1, -1\\4\nmid\tilde{m}_\xi}}
|\{w_\xi\in Z^{\nu_\xi}_\inv\,|\,
\ell_{w_\xi}^{1}(\nu_\xi)=\ell_{w_\xi}^{2}(\nu_\xi)=0\}|\\
&\thickspace\times
\prod_{\substack{\{\xi_1\neq\xi_2\}\in\langle\tsigma\rangle\setminus L\\
\xi_1^{\vee}=\xi_2\\ 2|\tilde{m}_{\xi_1}}} 
\delta_{\nu_{\xi_1}\nu_{\xi_2}}z_{\nu_{\xi_1}}
\prod_{\substack{\{\xi_1\neq\xi_2\}\in\langle\tsigma\rangle\setminus L\\
\xi_1^{\vee}=\xi_2\\ 2\nmid\tilde{m}_{\xi_1}}}
(-1)^{\ell(\nu_{\xi_1})_1}
\delta_{\nu_{\xi_1}\nu_{\xi_2}}z_{\nu_{\xi_1}}.
\end{split}
\end{equation*}
In the second factor, we have noted that if 
$w_{-1}\in Z^{\nu_{-1}}_{\ff-\inv}$, then
\begin{equation*}
\begin{split}
\ell_{w_{-1}}^{2}(\nu_{-1})_{0\,\mathrm{mod}\,4}
+ \frac{1}{2}\ell_{w_{-1}}^{3}(\nu_{-1})_1 &\equiv
\ell(\nu_{-1})_{0\,\mathrm{mod}\,4} + \frac{1}{2}\ell(\nu_{-1})_1\\
&\equiv \frac{|\nu_{-1}|}{2} + \ell(\nu_{-1})_0 \ \mathrm{mod}\ 2.
\end{split}
\end{equation*}
In the third factor, since $4|\tilde{m}_\xi$, both
$\Lambda_{w_\xi}^{1}(\nu_\xi)$ and $\Lambda_{w_\xi}^{2}(\nu_\xi)$
contribute to $\ell_w^{2}(\underline{\nu})_{0\,\mathrm{mod}\,4}$,
and nothing contributes to $\ell_w^{3}(\underline{\nu})_1$.
In the fourth factor, since $2|\tilde{m}_\xi$ but $4\nmid\tilde{m}_\xi$,
and $\ell_{w_\xi}^{1}(\nu_\xi)_0=\ell_{w_\xi}^{2}(\nu_\xi)_0=0$,
there is no contribution to the sign.
So in addition to \eqref{combglnspneqn} and \eqref{othercombglnoneqn}, 
we need the following fact:
\begin{equation}
\label{combunununeqn}
\sum_{(w,\epsilon)\in Z^{\nu}_{(p^{+},p^{-})-\inv}}
\negthickspace\negthickspace
(-1)^{\ell_w^{2}(\nu)_{0\,\mathrm{mod}\,4}+\frac{1}{2}\ell_w^{3}(\nu)_1}
=\sum_{\rho\vdash|\nu|} (-1)^{n(\rho)}
|\T_{(p^{+},p^{-})}(\rho')^{\psi}|\chi_\nu^{\rho}.
\end{equation}
This will be proved in \S5.
\subsection{The $U_n(\Fqs)/U_{n/2}(\mathbb{F}_{q^4})$ Case}
The final case to consider in this section is when $F:G\to G$
is a non-split Frobenius map for which $F(s)=-s$. 
Replacing $s$ by a $G^{F}$-conjugate
if necessary, we may assume that there is a form 
$\langle\cdot,\cdot\rangle$ on $V$ and an
involution $\theta':G\to G$ with the same properties as in \S3.3. Then
$\theta'F$ is the split Frobenius map induced by some $F_V:V\to V$
which respects $\langle\cdot,\cdot\rangle$ and interchanges 
$V^{+}$ and $V^{-}$. In particular, $n$ is even, and
$n^{+}=n^{-}=\frac{n}{2}$.
One has $G^{F}\cong U_n(\Fqs)$, $K^{F}\cong U_{n/2}(\mathbb{F}_{q^4})$.
The result is:
\begin{theorem} \label{ununthm}
For any $\underline{\rho}\in\widehat{\CalP}_{n}^{\tsigma}$,
\[ \langle \chi^{\underline{\rho}},
\Ind_{U_{n/2}(\mathbb{F}_{q^4})}^{U_n(\Fqs)}(1)\rangle=
\left\{
\begin{array}{cl}
{\displaystyle\prod_i(m_{2i}(\rho_1')+1)}, &\begin{array}{l}
\text{if $\rho_\xi=\rho_{\xi^{-1}}, \forall\xi$,}\\
2|m_{2i+1}(\rho_1'), \forall i,\\
\text{and $\rho_{-1}'$ is even}\end{array}\\ 
0, &\text{ otherwise.} 
\end{array} \right. \]
\end{theorem}
By \eqref{uneqn}, it is equivalent to say that for any
$\underline{\nu}\in\widehat{\CalP}_{n}^{\tsigma}$,
\begin{equation} \label{ununeqn}
\begin{split}
\langle B_{\underline{\nu}}, \Ind_{K^{F}}^{G^{F}}(1)\rangle&=
(\sum_{\substack{\rho_1\vdash|\nu_1|\\
2|m_{2i+1}(\rho_1')}}(-1)^{n(\rho_1)}
(\prod_i(m_{2i}(\rho_1')+1))\chi_{\nu_1}^{\rho_1})\\
&\qquad
\thickspace\thickspace\times((-1)^{\frac{|\nu_{-1}|}{2}}
\sum_{\substack{\rho_{-1}\vdash|\nu_{-1}|\\\text{$\rho_{-1}'$ even}}}
\chi_{\nu_{-1}}^{\rho_{-1}})\\
&\qquad\times
\prod_{\substack{\xi\in\langle\tsigma\rangle\setminus L\\
\xi^{\vee}=\xi\\ \xi\neq 1, -1\\ 4|\tilde{m}_\xi}}
((-1)^{|\nu_\xi|}\sum_{\rho_\xi\vdash|\nu_\xi|} \chi_{\nu_\xi}^{\rho_\xi})
\prod_{\substack{\xi\in\langle\tsigma\rangle\setminus L\\
\xi^{\vee}=\xi\\ \xi\neq 1, -1\\ 4\nmid\tilde{m}_\xi}}
(\sum_{\rho_\xi\vdash|\nu_\xi|} \chi_{\nu_\xi}^{\rho_\xi})\\
&\qquad\times
\prod_{\substack{\{\xi_1\neq\xi_2\}\in\langle\tsigma\rangle\setminus L\\
\xi_1^{\vee}=\xi_2\\ 2|\tilde{m}_{\xi_1}}} 
\delta_{\nu_{\xi_1}\nu_{\xi_2}}z_{\nu_{\xi_1}}
\prod_{\substack{\{\xi_1\neq\xi_2\}\in\langle\tsigma\rangle\setminus L\\
\xi_1^{\vee}=\xi_2\\ 2\nmid\tilde{m}_{\xi_1}}} 
(-1)^{|\nu_{\xi_1}|}\delta_{\nu_{\xi_1}\nu_{\xi_2}}z_{\nu_{\xi_1}}.
\end{split}
\end{equation}
For the signs here, see the comments after \eqref{unununeqn}.

\newpar
Let $(T,\lambda)$ be as in \S3.3. We have:
\begin{lemma}
The map $f\mapsto (w_f,\epsilon_f)$ induces a bijection
\[ T^{F}\!\setminus\!\Theta_T^{F}/K^{F}\isomto
W(T)^{F}_{\star-\inv}. \]
\end{lemma}
Now Lemma \ref{unununlemma} holds again here, with the same 
$Y^{\underline{\nu}}_{\inv}$
and the same proof. So arguing as in \S3.3, we get
\begin{equation*}
\langle B_{\underline{\nu}}, \Ind_{K^{F}}^{G^{F}}(1)\rangle
=\negthickspace\negthickspace\negthickspace\negthickspace
\sum_{(w,\epsilon)\in Y^{\underline{\nu}}_{\star-\inv}}
\negthickspace\negthickspace\negthickspace\negthickspace
(-1)^{l_w^{2}(\underline{\nu})_{0\,\mathrm{mod}\,4}
+\frac{1}{2}\ell_w^{3}(\underline{\nu})_1}.
\end{equation*}
The rest of the proof is also the same as in \S3.3, except that Equation
\eqref{combunununeqn} is replaced by:
\begin{equation} \label{combununeqn}
\begin{split}
\sum_{(w,\epsilon)\in Z^{\nu}_{\star-\inv}}\negthickspace\negthickspace
(-1)&^{l_w^{2}(\nu)_{0\,\mathrm{mod}\,4}+\frac{1}{2}\ell_w^{3}(\nu)_1}\\
&=\sum_{\substack{\rho\vdash|\nu|\\2|m_{2i+1}(\rho')}} 
(-1)^{n(\rho)}
(\prod_i (m_{2i}(\rho')+1))\chi_\nu^{\rho}.
\end{split}
\end{equation}
This will be proved in \S5.
\section{Cases where $G/K=GL_n/O_n$ or $GL_n/SO_n$}
\setcounter{theorem}{0}
In this section, we suppose that $V$ has a nondegenerate
symmetric form $\langle\cdot,\cdot\rangle$, and that
$\theta:G\to G$ is the involution defined by
\[ \langle\theta(g)v,v'\rangle=\langle v,g^{-1}v'\rangle,\ \forall g\in G,
v,v'\in V. \]
Since $G^{\theta}=O(V,\langle\cdot,\cdot\rangle)$ has two components,
$K$ can be either $O_n$ or $SO_n$. It suffices to solve the problem
for $O_n$:
\begin{lemma} \label{sonlemma}
Let $F:G\to G$ be a Frobenius morphism which commutes with $\theta$.
Let $\underline{\rho}\in\widehat{\CalP}_n^{\sigma}$ if $F$ is split,
$\underline{\rho}\in\widehat{\CalP}_n^{\tsigma}$ if $F$ is non-split.
If $F$ is split, let $\zeta\in L^{\sigma}$ be such that 
$\langle -1,\zeta\rangle^{\sigma}=-1$. 
If $F$ is non-split, let $\zeta\in L^{\tsigma}$
be such that $\langle -1,\zeta\rangle^{\tsigma}=-1$. Then
\[ \langle\chi^{\underline{\rho}},\Ind_{((G^{\theta})^{\circ})^{F}}^{G^{F}}(1)
\rangle = \langle\chi^{\underline{\rho}},
\Ind_{(G^{\theta})^{F}}^{G^{F}}(1)\rangle +
\langle\chi^{\zeta.\underline{\rho}},
\Ind_{(G^{\theta})^{F}}^{G^{F}}(1)\rangle. \]
\end{lemma}
\begin{proof}
Clearly
\[ \Ind_{((G^{\theta})^{\circ})^{F}}^{(G^{\theta})^{F}}(1) =
\mathrm{Res}_{(G^{\theta})^{F}}^{G^{F}}
(1+\langle\det(\cdot),\zeta^{-1}\rangle), \]
so
\[ \Ind_{((G^{\theta})^{\circ})^{F}}^{G^{F}}(1) =
(1+\langle\det(\cdot),\zeta^{-1}\rangle).\Ind_{(G^{\theta})^{F}}^{G^{F}}(1), \]
which proves the result.
\end{proof}
For the remainder of this section, we will write $K$ for $G^{\theta}\cong O_n$
and $K^{\circ}$ for $(G^{\theta})^{\circ}\cong SO_n$.

\newpar
Let $T$ be a maximal torus of $G$. 
For $f\in\Theta_T$, we have $w_f\in W(T)_\inv$. This time
\[ T\cap fKf^{-1} = \{t\in T\,|\, w_f(t)=t^{-1}\} \]
is not necessarily connected.
\begin{proposition} \label{glnonprop}
The map $f\mapsto w_f$ is a bijection
$T\!\setminus\!\Theta_T/K\to W(T)_\inv$. Moreover, if $w\in W(T)_\inv$,
the corresponding $T$--$K$ double coset breaks into two
$T$--$K^{\circ}$ double cosets if $w\in W(T)_{\ff-\inv}$, and is a
single $T$--$K^{\circ}$ double coset otherwise.
\end{proposition}
\begin{proof}
As with Propositions \ref{glnspnprop} and \ref{glnglnglnprop},
this is better known as a statement about $K$-orbits on the flag variety
(see \cite[\S6]{trapa}). It is easy to prove.
\end{proof}
\subsection{The $GL_n(\Fq)/O_n(\Fq)$ Case
($n$ odd)}
In this subsection, suppose that $n$ is odd and
let $F:G\to G$ be a split Frobenius map
which commutes with $\theta$. 
So $F$ is induced by a Frobenius map $F_V$ on $V$ which respects
$\langle\cdot,\cdot\rangle$, such that $\langle\cdot,\cdot\rangle$
has Witt index $\lfloor\frac{n}{2}\rfloor$ on $V^{F_V}$. One has
$G^{F}\cong GL_n(\Fq)$, $K^{F}\cong O_n(\Fq)$. 
Recall the definition of $d_\xi$ from \S1.2. The result is:
\begin{theorem} \label{glnonthm}
For any $\underline{\rho}\in\widehat{\CalP}_n^{\sigma}$,
\begin{equation*}
\langle \chi^{\underline{\rho}},
\Ind_{O_n(\Fq)}^{GL_n(\Fq)}(1)\rangle =
\left\{\begin{array}{cl}
{\displaystyle\frac{1}{2}
\prod_{\substack{\xi\in\langle\sigma\rangle\setminus L\\
d_\xi=1}} (\prod_i (m_i(\rho_\xi)+1))},
&\text{ if $d_\xi=-1\Rightarrow\rho_\xi'$ is even}\\
0, &\text{ otherwise.}
\end{array}\right.
\end{equation*}
\end{theorem}
By \eqref{glneqn}, it is equivalent to say that for any
$\underline{\nu}\in\widehat{\CalP}_n^{\sigma}$,
\begin{equation} \label{glnoneqn}
\begin{split}
\langle B_{\underline{\nu}}, \Ind_{K^{F}}^{G^{F}}(1)\rangle&=
-\frac{1}{2}\prod_{\substack{\xi\in\langle\sigma\rangle\setminus L\\
d_\xi=1}}((-1)^{|\nu_\xi|}\sum_{\rho_\xi\vdash|\nu_\xi|}
(\prod_i (m_i(\rho_\xi)+1))\chi_{\nu_\xi}^{\rho_\xi})\\
&\qquad\thickspace\thickspace\times
\prod_{\substack{\xi\in\langle\sigma\rangle\setminus L\\ d_\xi=-1}}
(\sum_{\substack{\rho_\xi\vdash|\nu_\xi|\\\rho_\xi' \text{ even}}}
\chi_{\nu_\xi}^{\rho_\xi}).
\end{split}
\end{equation}
Here the sign in \eqref{glneqn} has been distributed in an obvious way.

\newpar
Let $(T,\lambda)$ be a pair in the $G^{F}$-orbit corresponding to
$\underline{\nu}\in\widehat{\CalP}_{n}^{\sigma}$.
\begin{lemma} \label{glnontrivlemma}
The map $f\mapsto w_f$ induces a surjection
\[ T^{F}\!\setminus\!\Theta_T^{F}/K^{F}\twoheadrightarrow
W(T)^{F}_\inv. \]
\end{lemma}
\begin{proof}
By Proposition \ref{glnonprop} and Lemma \ref{glnontrivlemma},
we may identify $W(T)^{F}_\inv$
with $(T\!\setminus\!\Theta_T/K)^{F}$. Under this identification,
an involution is in the image of $T^{F}\!\setminus\!\Theta_T^{F}/K^{F}$
precisely when the corresponding $F$-stable $K$-orbit
on $T\!\setminus\!\Theta_T$ contains an $F$-fixed point (by
connectedness of $T$). Since $n$ is odd, this $K$-orbit is a single
$K^{\circ}$-orbit, so this is automatic.
\end{proof}

\newpar
\begin{lemma} \label{glnondoublelemma}
For $f\in\Theta_T^{F}$, there are $2^{\ell_{w_f}^{1}(\underline{\nu})-1}$
$T^{F}$--$K^{F}$ double cosets in $(TfK)^{F}$.
\end{lemma}
\begin{proof}
Since the image of the Lang map on $K$ is $K^{\circ}$, the number
of $T^{F}$--$K^{F}$ double cosets in $(TfK)^{F}$ is the same as the number of
orbits of $T\cap fKf^{-1}$ on $T\cap fK^{\circ}f^{-1}$ for the
action $t.t'=tt'F(t)^{-1}$. Using Lang's theorem, one sees that only
the ``disconnected part'' of $T\cap fKf^{-1}$ contributes. So
this is the same as the number of orbits of
\[ 
\prod_{(\xi,j)\in\Lambda_{w_f}^{1}(\underline{\nu})}\negthickspace
(\pm 1)^{m_\xi(\nu_\xi)_j}
\text{ on }
\left\{(\epsilon_{(\xi,j,i)}')\in
\negthickspace\negthickspace
\prod_{(\xi,j)\in\Lambda_{w_f}^{1}(\underline{\nu})}\negthickspace
(\pm 1)^{m_\xi(\nu_\xi)_j}\,\left|\,
\prod_{(\xi,j,i)}\epsilon_{(\xi,j,i)}'=1\right.\right\} \]
for the action $(\epsilon_{(\xi,j,i)}).(\epsilon_{(\xi,j,i)}')
=(\epsilon_{(\xi,j,i)}\epsilon_{(\xi,j,i-1)}\epsilon_{(\xi,j,i)}')$.
In other words, it is the index in the latter group of the subgroup
where $\prod_i\epsilon_{(\xi,j,i)}'=1$, for all 
$(\xi,j)\in\Lambda_{w_f}^{1}(\underline{\nu})$.
Now since $n$ is odd, $\Lambda_{w_f}^{1}(\underline{\nu})$ is non-empty.
The result follows.
\end{proof}

\newpar
Now let $X^{\underline{\nu}}_\inv=\{w\in W(T)^{F}_{\lambda,\inv}\,|\,
(\xi,j)\in \Lambda_{w}^{1}(\underline{\nu})\Rightarrow
\langle -1,\xi\rangle^{\sigma^{m_\xi(\nu_\xi)_j}}=1\}$.
\begin{lemma} \label{glnonlemma}
For $f\in\Theta_T^{F}$, $f\in\Theta_{T,\lambda}^{F}\Leftrightarrow
w_f\in X^{\underline{\nu}}_\inv$.
\end{lemma}
\begin{proof}
We first note that $\epsilon_{T,f}=1$, as may be seen directly
(using formulas for $\Fq$-rank such as those below) or deduced by the
method of \cite[Lemma 11.3]{symmfinite}. So
\begin{equation*}
\Theta_{T,\lambda}^{F}=\{f\in\Theta_T^{F}\,|\,
\lambda|_{\{t\in T^{F}|w_f(t)=t^{-1}\}}=1\}.
\end{equation*}
Thus it suffices to observe that
$\{t\in T^{F}\,|\,w_f(t)=t^{-1}\}$ is generated by
\[ \{tw_f(t)\,|\,t\in T^{F}\}\text{ and }
\prod_{(\xi,j)\in\Lambda_{w_f}^{1}(\underline{\nu})}(\pm 1). \]
\end{proof}
\begin{corollary}
The map $f\mapsto w_f$ induces a surjection
$T^{F}\!\setminus\!\Theta_{T,\lambda}^{F}/K^{F}\twoheadrightarrow 
X^{\underline{\nu}}_\inv$.
\end{corollary}

\newpar
Now the $\Fq$-rank of $T$ is $\ell(\underline{\nu})$, and that of
\begin{equation*}
\begin{split}
Z_G((T\cap fKf^{-1})^{\circ})
&=GL\left(\bigoplus_{\substack{(\xi,j,i)\\
(\xi,j)\in\Lambda_{w_f}^{1}(\underline{\nu})}} L_{(\xi,j,i)}\right)\\
&\thickspace\times\prod_{\substack{(\xi,j,i)\\
(\xi,j)\in\Lambda_{w_f}^{2}(\underline{\nu})}} 
\negthickspace\negthickspace GL(L_{(\xi,j,i)})
\times\prod_{\substack{(\xi,j,i)\\
(\xi,j)\in\Lambda_{w_f}^{3}(\underline{\nu})}} 
\negthickspace\negthickspace GL(L_{(\xi,j,i)})
\end{split}
\end{equation*}
is
\begin{equation*}
\ell(\underline{\nu})+\negthickspace\negthickspace
\sum_{\substack{(\xi,j,i)\\
(\xi,j)\in\Lambda_{w_f}^{1}(\underline{\nu})}} 
\negthickspace\negthickspace(m_\xi(\nu_\xi)_j-1)
\equiv \ell(\underline{\nu})+n+\ell_{w_f}^{1}(\underline{\nu})
\ \mathrm{mod}\ 2.
\end{equation*}
So Lusztig's formula gives
\begin{equation*}
\begin{split}
\langle B_{\underline{\nu}}, \Ind_{K^{F}}^{G^{F}}(1)\rangle
&=(-1)^{n}\negthickspace\negthickspace
\sum_{f\in T^{F}\setminus\Theta_{T,\lambda}^{F}/K^{F}}
(-1)^{\ell_{w_f}^{1}(\underline{\nu})}\\
&=-\sum_{w\in X^{\underline{\nu}}_\inv}
(-1)^{\ell_{w}^{1}(\underline{\nu})}
2^{\ell_{w}^{1}(\underline{\nu})-1}\\
&=-\frac{1}{2}\sum_{w\in X^{\underline{\nu}}_\inv}
(-2)^{\ell_{w}^{1}(\underline{\nu})}.
\end{split}
\end{equation*}

\newpar
Under the isomorphism $W(T)^{F}_{\lambda}\isomto
\prod_{\xi\in\langle\sigma\rangle\!\setminus\! L}Z^{\nu_\xi}$
of \S1.4, $X^{\underline{\nu}}_\inv$ corresponds to
\[ \{(w_\xi)\in\prod_{\xi\in\langle\sigma\rangle\setminus L}
Z^{\nu_\xi}_\inv\,|\, d_\xi=-1 \Rightarrow
\ell_{w_\xi}^{1}(\nu_\xi)_1=0 \}. \]
Hence
\begin{equation*}
\langle B_{\underline{\nu}}, \Ind_{K^{F}}^{G^{F}}(1)\rangle
=-\frac{1}{2}\prod_{\xi\in\langle\sigma\rangle\setminus L}
\left(\sum_{\substack{
w_\xi\in Z^{\nu_\xi}_\inv\\
d_\xi=-1 \Rightarrow
\ell_{w_\xi}^{1}(\nu_\xi)_1=0}}
\negthickspace\negthickspace
(-2)^{\ell_{w_\xi}^{1}(\nu_\xi)}\right).
\end{equation*}
So to prove \eqref{glnoneqn}, the combinatorial facts we need are:
\begin{gather}
\label{combglnoneqn}
\sum_{w\in Z^{\nu}_\inv}
(-2)^{\ell_{w}^{1}(\nu)}
=(-1)^{|\nu|}\sum_{\rho\vdash|\nu|}
(\prod_i (m_i(\rho)+1))\chi_\nu^{\rho},\text{ and}
\\
\label{combglnsoneqn}
\sum_{\substack{w\in Z^{\nu}_\inv\\\ell_w^{1}(\nu)_1=0}}
\negthickspace\negthickspace\negthickspace
(-2)^{\ell_w^{1}(\nu)}
=\sum_{\substack{\rho\vdash|\nu|\\\rho' \text{ even}}}
\chi_\nu^{\rho}. 
\end{gather}
These will be proved in \S5.
\subsection{The $GL_n(\Fq)/O_n^{\pm}(\Fq)$ Case ($n$ even)}
In this subsection, suppose that $n$ is even and
$F:G\to G$ is a split Frobenius map
which commutes with $\theta$.
So $F$ is induced by a Frobenius map $F_V$ on $V$ which respects
$\langle\cdot,\cdot\rangle$, and once again
$G^{F}\cong GL_n(\Fq)$.
The Witt index of $\langle\cdot,\cdot\rangle$ on $V^{F_V}$
is either $\frac{n}{2}$ or $\frac{n}{2}-1$, and accordingly either
$K^{F}\cong O_n^{+}(\Fq)$ or
$K^{F}\cong O_n^{-}(\Fq)$. Let $\epsilon$ be
$1$ in the first case and $-1$ in the second case.
\begin{theorem} \label{glnonpmthm}
For any $\underline{\rho}\in\widehat{\CalP}_n^{\sigma}$,
\begin{equation*}
\begin{split}
\langle \chi^{\underline{\rho}},
\Ind_{O_n^{\epsilon}(\Fq)}^{GL_n(\Fq)}(1)\rangle &=
\left\{\begin{array}{cl}
{\displaystyle\frac{1}{2}
\prod_{\substack{\xi\in\langle\sigma\rangle\setminus L\\
d_\xi=1}} (\prod_i (m_i(\rho_\xi)+1))},
&\text{ if $d_\xi=-1\Rightarrow\rho_\xi'$ is even}\\
0, &\text{ otherwise}
\end{array}\right.\\
&\thickspace+\left\{\begin{array}{cl}
\frac{1}{2}\epsilon,
&\text{ if all $\rho_\xi'$ are even}\\
0, &\text{ otherwise.}
\end{array}\right.
\end{split}
\end{equation*}
\end{theorem}
By \eqref{glneqn}, it is equivalent to say that for any
$\underline{\nu}\in\widehat{\CalP}_n^{\sigma}$,
\begin{equation} \label{glnonpmeqn}
\begin{split}
\langle B_{\underline{\nu}}, \Ind_{K^{F}}^{G^{F}}(1)\rangle&=
\frac{1}{2}\prod_{\substack{\xi\in\langle\sigma\rangle\setminus L\\
d_\xi=1}}((-1)^{|\nu_\xi|}\sum_{\rho_\xi\vdash|\nu_\xi|}
(\prod_i (m_i(\rho_\xi)+1))\chi_{\nu_\xi}^{\rho_\xi})\\
&\qquad\thickspace\thickspace\times
\prod_{\substack{\xi\in\langle\sigma\rangle\setminus L\\ d_\xi=-1}}
(\sum_{\substack{\rho_\xi\vdash|\nu_\xi|\\\rho_\xi' \text{ even}}}
\chi_{\nu_\xi}^{\rho_\xi})\\
&\thickspace+\frac{1}{2}\epsilon\prod_{\xi\in\langle\sigma\rangle\setminus L}
(\sum_{\substack{\rho_\xi\vdash|\nu_\xi|\\\rho_\xi' \text{ even}}}
\chi_{\nu_\xi}^{\rho_\xi}).
\end{split}
\end{equation}

\newpar
Let $(T,\lambda)$ be a pair in the $G^{F}$-orbit corresponding to
$\underline{\nu}\in\widehat{\CalP}_{n}^{\sigma}$.
\begin{lemma} \label{glnonpmtrivlemma}
The map $f\mapsto w_f$ induces a map
$T^{F}\!\setminus\!\Theta_T^{F}/K^{F}\rightarrow
W(T)^{F}_\inv$.
If $\epsilon_{\underline{\nu}}=\epsilon$, this map is surjective.
If $\epsilon_{\underline{\nu}}=-\epsilon$, the image is
$W(T)^{F}_\inv\setminus W(T)^{F}_{\ff-\inv}$.
\end{lemma}
\begin{proof}
As in the proof of Lemma \ref{glnontrivlemma},
an involution is in the image of $T^{F}\!\setminus\!\Theta_T^{F}/K^{F}$
precisely when the corresponding $F$-stable $K$-orbit
on $T\!\setminus\!\Theta_T$ contains an $F$-fixed point. 
If the involution has a fixed point, this $K$-orbit is a single
$K^{\circ}$-orbit, so this is automatic. Suppose $w\in W(T)^{F}_{\ff-\inv}$.
It is in the image of
$T^{F}\!\setminus\!\Theta_T^{F}/K^{F}$ precisely when there exists
a decomposition of $V$ into lines $\{L_i'\,|\,1\leq i\leq n\}$ such that
\begin{enumerate}
\item  $(L_i')^{\perp}=\oplus_{i'\neq w(i)} L_{i'}'$, and
\item  $F_V(L_i')=L_{w_{\underline{\nu}}(i)}'$
\end{enumerate}
(here we have identified $W(T)$ with $S_n$ in some way). We must prove
that this happens if and only if 
$\mathrm{sign}(w_{\underline{\nu}})=\epsilon$.
Since $\epsilon$ is multiplicative with respect to $F_V$-stable
orthogonal direct sums, we may assume that $\langle w\rangle\times
\langle w_{\underline{\nu}}\rangle$ acts transitively on $\{1,\cdots,n\}$.
So if $\mathrm{sign}(w_{\underline{\nu}})=1$, then 
$\langle w_{\underline{\nu}}\rangle$ has two orbits
on $\{1,\cdots,n\}$ which $w$ interchanges; in this case the existence
of $\{L_i'\}$ as above is clearly equivalent to the existence of a
decomposition $V=V_1\oplus V_2$ into $F_V$-stable Lagrangian subspaces,
which indeed happens if and only if $\epsilon=1$. On the other hand,
if $\mathrm{sign}(w_{\underline{\nu}})=-1$, then 
$w_{\underline{\nu}}$ is an $n$-cycle and
$w=w_{\underline{\nu}}^{n/2}$; in this case the existence
of $\{L_i'\}$ as above is easily seen to be equivalent to the existence of a
Lagrangian subspace $V_1$ of $V$ such that $\dim V_1\cap F_V(V_1)
=\frac{n}{2}-1$, which indeed happens if and only if $\epsilon=-1$. 
\end{proof}

\newpar
\begin{lemma} \label{glnonpmdoublelemma}
For $f\in\Theta_T^{F}$, the number of $T^{F}$--$K^{F}$ double cosets
in $(TfK)^{F}$ is
\[ \left\{\begin{array}{cl}
1, &
\text{if $w_f\in W(T)^{F}_{\ff-\inv}$}\\
2^{\ell_{w_f}^{1}(\underline{\nu})-1}, 
&\text{ otherwise.}\end{array}\right.\]
\end{lemma}
\begin{proof}
The method of proof of Lemma \ref{glnondoublelemma} applies again here.
\end{proof}

\newpar
Now define $X^{\underline{\nu}}_\inv$ as in \S4.1. 
Let $X^{\underline{\nu}}_{\ff-\inv}=X^{\underline{\nu}}_\inv
\cap W(T)^{F}_{\ff-\inv}$.
\begin{lemma} \label{glnonpmlemma}
For $f\in\Theta_T^{F}$, $f\in\Theta_{T,\lambda}^{F}\Leftrightarrow
w_f\in X^{\underline{\nu}}_\inv$.
\end{lemma}
\begin{proof}
The proof is exactly the same as that of Lemma \ref{glnonlemma}.
\end{proof}
\begin{corollary}
\begin{enumerate}
\item  The map $f\mapsto w_f$ induces a map
$T^{F}\!\setminus\!\Theta_{T,\lambda}^{F}/K^{F}\to X^{\underline{\nu}}_\inv$.
If $\epsilon_{\underline{\nu}}=\epsilon$, this map is
surjective; if $\epsilon_{\underline{\nu}}=-\epsilon$, its image is
$X^{\underline{\nu}}_\inv\setminus X^{\underline{\nu}}_{\ff-\inv}$.
\item  If $w$ is in the image of 
$T^{F}\!\setminus\!\Theta_{T,\lambda}^{F}/K^{F}\to X^{\underline{\nu}}_\inv$,
there are
\[ \left\{\begin{array}{cl}
1, &\text{ if $w\in X^{\underline{\nu}}_{\ff-\inv}$}\\ 
2^{\ell_{w_f}^{1}(\underline{\nu})-1}, 
&\text{ otherwise}\end{array}\right.\]
$T^{F}$--$K^{F}$ double cosets in the preimage of $w$.
\end{enumerate}
\end{corollary}
\begin{proof}
This follows by combining 
Proposition \ref{glnonprop}, Lemma \ref{glnonpmtrivlemma},
Lemma \ref{glnonpmdoublelemma} and Lemma \ref{glnonpmlemma}.
\end{proof}

\newpar
Now as in \S4.1, the $\Fq$-rank of $T$ is $\ell(\underline{\nu})$, and that of
$Z_G((T\cap fKf^{-1})^{\circ})$
is congruent to
$\ell(\underline{\nu})+n+\ell_{w_f}^{1}(\underline{\nu})\ \mathrm{mod}\ 2$.
So Lusztig's formula gives
\begin{equation*}
\begin{split}
\langle B_{\underline{\nu}}, \Ind_{K^{F}}^{G^{F}}(1)\rangle
&=\sum_{f\in T^{F}\setminus\Theta_{T,\lambda}^{F}/K^{F}}
(-1)^{\ell_{w_f}^{1}(\underline{\nu})}\\
&=\sum_{w\in X^{\underline{\nu}}_\inv}
(-1)^{\ell_{w}^{1}(\underline{\nu})}
\left\{\begin{array}{cl}
1, &\text{ if $w\in X^{\underline{\nu}}_{\ff-\inv}$}\\
2^{\ell_{w}^{1}(\underline{\nu})-1},
&\text{ otherwise}\end{array}\right.\\
&\qquad\thickspace\thickspace\thickspace
-\frac{1}{2}(1-\epsilon\epsilon_{\underline{\nu}})
|X^{\underline{\nu}}_{\ff-\inv}|\\
&=\frac{1}{2}\sum_{w\in X^{\underline{\nu}}_\inv}
(-2)^{\ell_{w}^{1}(\underline{\nu})}
+\frac{1}{2}\epsilon_{\underline{\nu}}
\epsilon|X^{\underline{\nu}}_{\ff-\inv}|.
\end{split}
\end{equation*}

\newpar
We transform this expression as in the previous subsection to obtain:
\begin{equation*}
\begin{split}
\langle B_{\underline{\nu}}, \Ind_{K^{F}}^{G^{F}}(1)\rangle
&=\frac{1}{2}\prod_{\xi\in\langle\sigma\rangle\setminus L}
\left(\sum_{\substack{
w_\xi\in Z^{\nu_\xi}_\inv\\
d_\xi=-1 \Rightarrow \ell_{w_\xi}^{1}(\nu_\xi)_1=0}}
\negthickspace\negthickspace
(-2)^{\ell_{w_\xi}^{1}(\nu_\xi)}\right)\\
&\thickspace+\frac{1}{2}\epsilon\prod_{\xi\in\langle\sigma\rangle\setminus L}
\epsilon_{\nu_\xi}|Z^{\nu_\xi}_{\ff-\inv}|.
\end{split}
\end{equation*}
In those factors of the second term for which $m_\xi$ is even,
we have used the fact that if there
exists a fixed-point free involution in $Z^{\nu_\xi}$,
then $\ell(\nu_\xi)_1$ is even, so 
$(-1)^{\ell(\nu_\xi)}=\epsilon_{\nu_\xi}$.
So along with \eqref{combglnoneqn} and \eqref{combglnsoneqn}, we need
\eqref{combglnspneqn} multiplied on both sides by $\epsilon_\nu$.
\subsection{The $U_n(\Fqs)/O_n(\Fq)$ Case ($n$ odd)}
In this subsection, suppose that $n$ is odd and
take $F:G\to G$ to be a non-split Frobenius map
commuting with $\theta$. 
Thus $\theta F$ is a split Frobenius map induced by $F_V$ as in \S4.1.
One has
$G^{F}\cong U_n(\Fqs)$, $K^{F}\cong O_n(\Fq)$. The result is:
\begin{theorem} \label{unonthm}
For any $\underline{\rho}\in\widehat{\CalP}_{n}^{\tsigma}$,
\begin{equation*}
\begin{split}
\langle \chi^{\underline{\rho}},
\Ind_{O_n(\Fq)}^{U_n(\Fqs)}(1)\rangle&=
\left\{
\begin{array}{rl}
{\displaystyle\frac{1}{2}
\prod_{\substack{\xi\in\langle\tsigma\rangle\setminus L\\
\tilde{d}_\xi=1\\2\nmid\tilde{m}_\xi}}}
&\negthickspace\negthickspace\negthickspace
{\displaystyle(\prod_i(m_{2i+1}(\rho_\xi)+1))
\prod_{\substack{\xi\in\langle\tsigma\rangle\setminus L\\
\tilde{d}_\xi=-1\\2\nmid\tilde{m}_\xi}}
(\prod_i(m_{2i}(\rho_\xi)+1))}\\
&\thickspace{\displaystyle\times
\prod_{\substack{\xi\in\langle\tsigma\rangle\setminus L\\
\tilde{d}_\xi=1\\2|\tilde{m}_\xi}}
(\prod_i(m_{i}(\rho_\xi)+1))},\\
&\begin{array}{l}
\text{if $\tilde{d}_\xi=1$,
$2\nmid \tilde{m}_\xi \Rightarrow 2|m_{2i}(\rho_\xi),\ \forall i$,}\\
\text{$\tilde{d}_\xi=-1$,
$2\nmid \tilde{m}_\xi \Rightarrow 2|m_{2i+1}(\rho_\xi),\ \forall i$,}\\
\text{and $\tilde{d}_\xi=-1$,
$2|\tilde{m}_\xi \Rightarrow \rho_\xi'$ is even}
\end{array}\\ 
0, &\text{ otherwise.} 
\end{array} \right.
\end{split}
\end{equation*}
\end{theorem}
By \eqref{uneqn}, it is equivalent to say that for any
$\underline{\nu}\in\widehat{\CalP}_{n}^{\tsigma}$,
\begin{equation} \label{unoneqn}
\begin{split}
\langle B_{\underline{\nu}}, \Ind_{K^{F}}^{G^{F}}(1)\rangle&=
\frac{1}{2}(-1)^{\lfloor\frac{n}{2}\rfloor}
\negthickspace\negthickspace
\prod_{\substack{\xi\in\langle\tsigma\rangle\setminus L\\
\tilde{d}_\xi=1\\2\nmid\tilde{m}_\xi}}
\left(\sum_{\substack{\rho_\xi\vdash|\nu_\xi|\\2|m_{2i}(\rho_\xi)}}
(-1)^{n(\rho_\xi')}
(\prod_i(m_{2i+1}(\rho_\xi)+1))\chi_{\nu_\xi}^{\rho_\xi}\right)\\
&\thickspace\thickspace\times
\prod_{\substack{\xi\in\langle\tsigma\rangle\setminus L\\
\tilde{d}_\xi=-1\\2\nmid\tilde{m}_\xi}}
\left(\sum_{\substack{\rho_\xi\vdash|\nu_\xi|\\2|m_{2i+1}(\rho_\xi)}}
(-1)^{n(\rho_\xi')}
(\prod_i(m_{2i}(\rho_\xi)+1))\chi_{\nu_\xi}^{\rho_\xi}\right)\\
&\thickspace\thickspace\times
\prod_{\substack{\xi\in\langle\tsigma\rangle\setminus L\\
\tilde{d}_\xi=1\\2|\tilde{m}_\xi}}
(-1)^{|\nu_\xi|}
\left(\sum_{\rho_\xi\vdash|\nu_\xi|}
(\prod_i(m_{i}(\rho_\xi)+1))\chi_{\nu_\xi}^{\rho_\xi}\right)\\
&\thickspace\thickspace\times
\prod_{\substack{\xi\in\langle\tsigma\rangle\setminus L\\
\tilde{d}_\xi=-1\\2|\tilde{m}_\xi}}
\left(\sum_{\substack{\rho_\xi\vdash|\nu_\xi|\\\rho_\xi'\text{ even}}}
\chi_{\nu_\xi}^{\rho_\xi}\right).
\end{split}
\end{equation}
Here the sign
\[ (-1)^{\lceil\frac{n}{2}\rceil+
\sum_{\xi\in\langle\tsigma\rangle\setminus L}
\tilde{m}_\xi n(\rho_\xi')+|\rho_\xi|}
= (-1)^{\lfloor\frac{n}{2}\rfloor+
\sum_{\xi\in\langle\tsigma\rangle\setminus L}
\tilde{m}_\xi(n(\rho_\xi')+|\rho_\xi|)+|\rho_\xi|} \]
of \eqref{uneqn} has been distributed in an obvious way.

\newpar
The proof of these statements is mostly very similar to that
of \eqref{glnoneqn}. Let $(T,\lambda)$ be a pair in the 
$G^{F}$-orbit corresponding to
$\underline{\nu}\in\widehat{\CalP}_{n}^{\tsigma}$.
\begin{lemma} \label{unontrivlemma}
The map $f\mapsto w_f$ induces a surjection
\[ T^{F}\!\setminus\!\Theta_T^{F}/K^{F}\twoheadrightarrow
W(T)^{F}_\inv. \]
\end{lemma}
\begin{proof}
This is deduced in the same way as Lemma \ref{glnontrivlemma}.
\end{proof}

\newpar
\begin{lemma} \label{unondoublelemma}
For $f\in\Theta_T^{F}$, there are $2^{\ell_{w_f}^{1}(\underline{\nu})-1}$
$T^{F}$--$K^{F}$ double cosets
in $(TfK)^{F}$.
\end{lemma}
\begin{proof}
The proof is the same as that of Lemma \ref{glnondoublelemma}.
\end{proof}

\newpar
Now define $X^{\underline{\nu}}_\inv$ as in \S4.1, but with
$\tsigma$ instead of $\sigma$ and $\tilde{m}_\xi$ instead of $m_\xi$.
\begin{lemma} \label{unonlemma}
For $f\in\Theta_T^{F}$, $f\in\Theta_{T,\lambda}^{F}\Leftrightarrow
w_f\in X^{\underline{\nu}}_\inv$.
\end{lemma}
\begin{proof}
We first prove that $\epsilon_{T,f}=1$. The $\Fq$-rank of
\begin{equation*}
\begin{split}
Z_G((T\cap fKf^{-1})^{\circ})
&=GL\left(\bigoplus_{\substack{(\xi,j,i)\\
(\xi,j)\in\Lambda_{w_f}^{1}(\underline{\nu})}}L_{(\xi,j,i)}\right)\\
&\qquad\thickspace\times
\negthickspace\negthickspace\negthickspace\negthickspace
\prod_{\substack{(\xi,j,i)\\
(\xi,j)\in\Lambda_{w_f}^{2}(\underline{\nu})}} 
\negthickspace\negthickspace\negthickspace\negthickspace
GL(L_{(\xi,j,i)})\thickspace
\times\negthickspace\negthickspace
\prod_{\substack{(\xi,j,i)\\
(\xi,j)\in\Lambda_{w_f}^{3}(\underline{\nu})}} 
GL(L_{(\xi,j,i)})
\end{split}
\end{equation*}
is
\begin{equation*}
\lfloor\frac{1}{2}\sum_{(\xi,j)\in\Lambda_{w_f}^{1}(\underline{\nu})}
\tilde{m}_\xi(\nu_\xi)_j\rfloor
+\ell_{w_f}^{2}(\underline{\nu})
+\ell_{w_f}^{3}(\underline{\nu})_0.
\end{equation*}
If $t\in(T\cap fKf^{-1})^{F}$ has eigenvalue $\alpha_{(\xi,j)}$
on $L_{(\xi,j,0)}$, then
\begin{equation*}
\begin{split}
Z_G^{\circ}(t)\cap Z_G((T&\cap fKf^{-1})^{\circ})\\
&=GL\left(\bigoplus_{\substack{(\xi,j,i)\\
(\xi,j)\in\Lambda_{w_f}^{1}(\underline{\nu})\\
\alpha_{(\xi,j)}=1}}L_{(\xi,j,i)}\right)
\times GL\left(\bigoplus_{\substack{(\xi,j,i)\\
(\xi,j)\in\Lambda_{w_f}^{1}(\underline{\nu})\\
\alpha_{(\xi,j)}=-1}}L_{(\xi,j,i)}\right) \\
&\qquad\thickspace\times
\negthickspace\negthickspace\negthickspace\negthickspace
\prod_{\substack{(\xi,j,i)\\
(\xi,j)\in\Lambda_{w_f}^{2}(\underline{\nu})}}
\negthickspace\negthickspace\negthickspace\negthickspace
GL(L_{(\xi,j,i)})\thickspace
\times\negthickspace\negthickspace
\prod_{\substack{(\xi,j,i)\\
(\xi,j)\in\Lambda_{w_f}^{3}(\underline{\nu})}} 
GL(L_{(\xi,j,i)})
\end{split}
\end{equation*}
has $\Fq$-rank which differs from that of $Z_G((T\cap fKf^{-1})^{\circ})$ by
\[ \lfloor\frac{1}{2}\sum_{\substack{
(\xi,j)\in\Lambda_{w_f}^{1}(\underline{\nu})\\
\alpha_{(\xi,j)}=1}}\tilde{m}_\xi(\nu_\xi)_j\rfloor
+ \lfloor\frac{1}{2}\sum_{\substack{
(\xi,j)\in\Lambda_{w_f}^{1}(\underline{\nu})\\
\alpha_{(\xi,j)}=-1}}\tilde{m}_\xi(\nu_\xi)_j\rfloor
- \lfloor\frac{1}{2}\sum_{(\xi,j)\in\Lambda_{w_f}^{1}(\underline{\nu})}
\tilde{m}_\xi(\nu_\xi)_j\rfloor. \]
Since $n$ is odd,
$\sum_{(\xi,j)\in\Lambda_{w_f}^{1}(\underline{\nu})}
\tilde{m}_\xi(\nu_\xi)_j$
is odd, so this difference is zero. Thus $\epsilon_{T,f}=1$. The rest
of the proof follows that of Lemma \ref{glnonlemma}.
\end{proof}
\begin{corollary}
The map $f\mapsto w_f$ induces a surjection
$T^{F}\!\setminus\!\Theta_{T,\lambda}^{F}/K^{F}\twoheadrightarrow 
X^{\underline{\nu}}_\inv$.
\end{corollary}

\newpar
Now the $\Fq$-rank of $T$ is $\ell(\underline{\nu})_0$,
and that of $Z_G((T\cap fKf^{-1})^{\circ})$ is given above, 
whence
\begin{equation*}
\begin{split}
\text{$\Fq$-rank}(T)&+\text{$\Fq$-rank}(Z_G((T\cap fKf^{-1})^{\circ}))\\
&\equiv \lfloor\frac{1}{2}\sum_{
(\xi,j)\in\Lambda_{w_f}^{1}(\underline{\nu})}
\tilde{m}_\xi(\nu_\xi)_j\rfloor
+\ell_{w_f}^{1}(\underline{\nu})_0\\
&\equiv \lfloor\frac{n}{2}\rfloor
+\ell_{w_f}^{1}(\underline{\nu})_0+
\ell_{w_f}^{2}(\underline{\nu})_{2\,\mathrm{mod}\,4}
+\frac{1}{2}\ell_{w_f}^{3}(\underline{\nu})_1\ \mathrm{mod}\ 2.
\end{split}
\end{equation*}
So Lusztig's formula gives (compare \S4.1):
\begin{equation*}
\langle B_{\underline{\nu}}, \Ind_{K^{F}}^{G^{F}}(1)\rangle
=\frac{1}{2}(-1)^{\lfloor\frac{n}{2}\rfloor}
\sum_{w\in X^{\underline{\nu}}_\inv}
(-1)^{\ell_{w}^{1}(\underline{\nu})_0+
\ell_{w}^{2}(\underline{\nu})_{2\,\mathrm{mod}\,4}
+\frac{1}{2}\ell_{w}^{3}(\underline{\nu})_1}
2^{\ell_w^{1}(\underline{\nu})}.
\end{equation*}

\newpar
As in \S4.1, we transform this expression to get
\begin{equation*}
\begin{split}
\langle B_{\underline{\nu}}, \Ind_{K^{F}}^{G^{F}}(1)\rangle
&=\frac{1}{2}(-1)^{\lfloor\frac{n}{2}\rfloor}\\
\times\prod_{\substack{\xi\in\langle\tsigma\rangle\setminus L\\
2\nmid\tilde{m}_\xi}}
&\left(\sum_{\substack{
w_\xi\in Z^{\nu_\xi}_\inv\\
\tilde{d}_\xi=-1 \Rightarrow 
\ell_{w_\xi}^{1}(\nu_\xi)_1=0}}
\negthickspace\negthickspace\negthickspace\negthickspace\negthickspace
(-1)^{\ell_{w_\xi}^{1}(\nu_\xi)_0+
\ell_{w_\xi}^{2}(\nu_\xi)_{2\,\mathrm{mod}\,4}
+\frac{1}{2}\ell_{w_\xi}^{3}(\nu_\xi)_1}
2^{\ell_{w_\xi}^{1}(\nu_\xi)}\right)\\
&\thickspace\times
\prod_{\substack{\xi\in\langle\tsigma\rangle\setminus L\\
2|\tilde{m}_\xi}}
\left(\sum_{\substack{
w_\xi\in Z^{\nu_\xi}_\inv\\
\tilde{d}_\xi=-1 \Rightarrow 
\ell_{w_\xi}^{1}(\nu_\xi)_1=0}}
\negthickspace\negthickspace
(-2)^{\ell_{w_\xi}^{1}(\nu_\xi)}\right).
\end{split}
\end{equation*}
In the second factor, we hagve noted that if $\tilde{m}_\xi$
is even, no element of $\Lambda_{w_\xi}^{2}(\nu_\xi)$
or $\Lambda_{w_\xi}^{3}(\nu_\xi)$ can contribute to
$\ell_{w}^{2}(\underline{\nu})_{2\,\mathrm{mod}\,4}
+\frac{1}{2}\ell_{w}^{3}(\underline{\nu})_1$, and every element of
$\Lambda_{w_\xi}^{1}(\nu_\xi)$ contributes to 
$\ell_{w}^{1}(\underline{\nu})_0$.
So in addition to \eqref{combglnoneqn} and \eqref{combglnsoneqn}, the
combinatorial facts we need are:
\begin{equation}
\label{combunoneqn}
\begin{split}
\sum_{w\in Z^{\nu}_\inv}&
\negthickspace\negthickspace
(-1)^{\ell_{w}^{1}(\nu)_0+\ell_{w}^{2}(\nu)_{2\,\mathrm{mod}\,4}
+\frac{1}{2}\ell_{w}^{3}(\nu)_1}
2^{\ell_{w}^{1}(\nu)}\\
&\qquad\qquad=\sum_{\substack{\rho\vdash|\nu|\\2|m_{2i}(\rho)}}
(-1)^{n(\rho')}
(\prod_i (m_{2i+1}(\rho)+1))\chi_\nu^{\rho},
\end{split}
\end{equation}
and
\begin{equation}
\label{combunsoneqn}
\begin{split}
\sum_{\substack{w\in Z^{\nu}_\inv\\
\ell_w^{1}(\nu)_1=0}}&
\negthickspace\negthickspace
(-1)^{\ell_{w}^{1}(\nu)+\ell_{w}^{2}(\nu)_{2\,\mathrm{mod}\,4}
+\frac{1}{2}\ell_{w}^{3}(\nu)_1}
2^{\ell_{w}^{1}(\nu)}\\
&\qquad\qquad=\sum_{\substack{\rho\vdash|\nu|\\2|m_{2i+1}(\rho)}}
(-1)^{n(\rho')}
(\prod_i (m_{2i}(\rho)+1))
\chi_\nu^{\rho}.
\end{split} 
\end{equation}
These will be proved in \S5.
\subsection{The $U_n(\Fqs)/O_n^{\pm}(\Fq)$ Case ($n$ even)}
Finally, we suppose that $n$ is even and
$F:G\to G$ is non-split.
So $\theta F$ is induced by $F_V$ as in \S4.2, and we have the same dichotomy
as to the Witt index of $\langle\cdot,\cdot\rangle$ on $V^{F_V}$.
Define $\epsilon\in\{\pm 1\}$ as in \S4.2, so that
$G^{F}\cong U_n(\Fqs)$,
$K^{F}\cong O_n^{\epsilon}(\Fq)$.
\begin{theorem} \label{unonpmthm}
For any $\underline{\rho}\in\widehat{\CalP}_{n}^{\tsigma}$,
\begin{equation*}
\begin{split}
\langle \chi^{\underline{\rho}},
\Ind_{O_n^{\epsilon}(\Fq)}^{U_n(\Fqs)}(1)\rangle&=
\left\{
\begin{array}{rl}
{\displaystyle\frac{1}{2}
\prod_{\substack{\xi\in\langle\tsigma\rangle\setminus L\\
\tilde{d}_\xi=1\\2\nmid\tilde{m}_\xi}}}
&\negthickspace\negthickspace\negthickspace
{\displaystyle(\prod_i(m_{2i}(\rho_\xi)+1))
\prod_{\substack{\xi\in\langle\tsigma\rangle\setminus L\\
\tilde{d}_\xi=-1\\2\nmid\tilde{m}_\xi}}
(\prod_i(m_{2i+1}(\rho_\xi)+1))}\\
&\thickspace{\displaystyle\times
\prod_{\substack{\xi\in\langle\tsigma\rangle\setminus L\\
\tilde{d}_\xi=1\\2|\tilde{m}_\xi}}
(\prod_i(m_{i}(\rho_\xi)+1)),}\\
&\begin{array}{l}
\text{if $\tilde{d}_\xi=1$,
$2\nmid \tilde{m}_\xi \Rightarrow 2|m_{2i+1}(\rho_\xi), \forall i$,}\\
\text{$\tilde{d}_\xi=-1$,
$2\nmid \tilde{m}_\xi \Rightarrow 2|m_{2i}(\rho_\xi), \forall i$,}\\
\text{and $\tilde{d}_\xi=-1$,
$2\,|\,\tilde{m}_\xi \Rightarrow \rho_\xi'$ is even}
\end{array}\\ 
0, &\text{ otherwise} 
\end{array} \right.\\
&\thickspace+\left\{\begin{array}{cl}
\frac{1}{2}\epsilon, &\text{ if all $\rho_\xi'$ are even}\\
0, &\text{ otherwise.}
\end{array}\right.
\end{split}
\end{equation*}
\end{theorem}
By \eqref{uneqn}, it is equivalent to say that for any
$\underline{\nu}\in\widehat{\CalP}_{n}^{\tsigma}$,
\begin{equation} \label{unonpmeqn}
\begin{split}
\langle B_{\underline{\nu}}, \Ind_{K^{F}}^{G^{F}}(1)\rangle&=
\frac{1}{2}(-1)^{\frac{n}{2}}\negthickspace
\prod_{\substack{\xi\in\langle\tsigma\rangle\setminus L\\
\tilde{d}_\xi=1\\2\nmid\tilde{m}_\xi}}
\left(\sum_{\substack{\rho_\xi\vdash|\nu_\xi|\\2|m_{2i+1}(\rho_\xi)}}
(-1)^{n(\rho_\xi')}
(\prod_i(m_{2i}(\rho_\xi)+1))\chi_{\nu_\xi}^{\rho_\xi}\right)\\
&\thickspace\thickspace\times
\prod_{\substack{\xi\in\langle\tsigma\rangle\setminus L\\
\tilde{d}_\xi=-1\\2\nmid\tilde{m}_\xi}}
\left(\sum_{\substack{\rho_\xi\vdash|\nu_\xi|\\2|m_{2i}(\rho_\xi)}}
(-1)^{n(\rho_\xi')}
(\prod_i(m_{2i+1}(\rho_\xi)+1))\chi_{\nu_\xi}^{\rho_\xi}\right)\\
&\thickspace\thickspace\times
\prod_{\substack{\xi\in\langle\tsigma\rangle\setminus L\\
\tilde{d}_\xi=1\\2|\tilde{m}_\xi}}
(-1)^{|\nu_\xi|}
\left(\sum_{\rho_\xi\vdash|\nu_\xi|}
(\prod_i(m_{i}(\rho_\xi)+1))\chi_{\nu_\xi}^{\rho_\xi}\right)\\
&\thickspace\thickspace\times
\prod_{\substack{\xi\in\langle\tsigma\rangle\setminus L\\
\tilde{d}_\xi=-1\\2|\tilde{m}_\xi}}
\left(\sum_{\substack{\rho_\xi\vdash|\nu_\xi|\\\rho_\xi'\text{ even}}}
\chi_{\nu_\xi}^{\rho_\xi}\right)\\
&\thickspace+
\frac{1}{2}\epsilon\prod_{\substack{\xi\in\langle\tsigma\rangle\setminus L
\\2\nmid\tilde{m}_\xi}}
(-1)^{\frac{|\nu_\xi|}{2}}
\left(\sum_{\substack{\rho_\xi\vdash|\nu_\xi|\\\rho_\xi'\text{ even}}}
\chi_{\nu_\xi}^{\rho_\xi}\right)
\prod_{\substack{\xi\in\langle\tsigma\rangle\setminus L
\\2|\tilde{m}_\xi}}
\left(\sum_{\substack{\rho_\xi\vdash|\nu_\xi|\\\rho_\xi'\text{ even}}}
\chi_{\nu_\xi}^{\rho_\xi}\right).
\end{split}
\end{equation}
Here the sign
\[ (-1)^{\frac{n}{2}+
\sum_{\xi\in\langle\tsigma\rangle\setminus L}
\tilde{m}_\xi n(\rho_\xi')+|\rho_\xi|}
= (-1)^{\frac{n}{2}+
\sum_{\xi\in\langle\tsigma\rangle\setminus L}
\tilde{m}_\xi(n(\rho_\xi')+|\rho_\xi|)+|\rho_\xi|} \]
of \eqref{uneqn} has been distributed in an obvious way
(in the second term it has been rewritten as
$\prod_{\xi\in\langle\tsigma\rangle\setminus L} 
(-1)^{\tilde{m}_\xi\frac{|\nu_\xi|}{2}}$ since $|\rho_\xi|$ and
$n(\rho_\xi')$ are even if $\rho_\xi'$ is even).

\newpar
The proof of this is similar to that of \eqref{glnonpmeqn}.
Let $(T,\lambda)$ be as in \S4.3.
Lemma \ref{glnonpmtrivlemma} must be modified as follows:
\begin{lemma} \label{unonpmtrivlemma}
The map $f\mapsto w_f$ induces a map
\[ T^{F}\!\setminus\!\Theta_T^{F}/K^{F}\to
W(T)^{F}_\inv. \]
If $\epsilon_{\underline{\nu}}=(-1)^{\frac{n}{2}}\epsilon$, 
this map is surjective.
If $\epsilon_{\underline{\nu}}=-(-1)^{\frac{n}{2}}\epsilon$, 
the image is $W(T)^{F}_\inv\setminus W(T)^{F}_{\ff-\inv}$.
\end{lemma}
\begin{proof}
The proof is mostly the same as that of Lemma \ref{glnonpmtrivlemma}.
Note that if $f\in\Theta_T^{F}$,
\begin{equation*}
\begin{split}
F_V(f^{-1}L_{(\xi,j,i)})&=\theta(f)^{-1}\bigcap_{\substack{
(\xi',j',i')\\\neq w_{\underline{\nu}}(\xi,j,i)}} 
L_{(\xi',j',i')}^{\perp}\\
&=\left(\bigoplus_{\substack{
(\xi',j',i')\\\neq w_{\underline{\nu}}(\xi,j,i)}} 
f^{-1}L_{(\xi',j',i')}\right)^{\perp}\\
&=f^{-1}L_{w_f^{-1}w_{\underline{\nu}}(\xi,j,i)}.
\end{split}
\end{equation*}
So in (2) of the proof of Lemma \ref{glnonpmtrivlemma}, 
$w_{\underline{\nu}}$
should be replaced by $w^{-1}w_{\underline{\nu}}$.
Since $\mathrm{sign}(w)=(-1)^{\frac{n}{2}}$ if 
$w\in W(T)^{F}_{\ff-\inv}$, we get the result.
\end{proof}

\newpar
\begin{lemma} \label{unonpmdoublelemma}
For $f\in\Theta_T^{F}$, the number of $T^{F}$--$K^{F}$ double cosets
in $(TfK)^{F}$ is
\[ \left\{\begin{array}{cl}
1, &
\text{if $w_f\in W(T)^{F}_{\ff-\inv}$}\\ 
2^{\ell_{w_f}^{1}(\underline{\nu})-1}, 
&\text{ otherwise.}\end{array}\right.\]
\end{lemma}
\begin{proof}
The proof of Lemma \ref{glnondoublelemma} applies again here.
\end{proof}

\newpar
Now we diverge somewhat from \S4.2, and define
\[ X^{\underline{\nu}}_\inv=
\{w\in W(T)_{\lambda,\inv}^{F}\,|\,(\xi,j)\in
\Lambda_{w}^{1}(\underline{\nu})\Rightarrow
\langle -1, \xi\rangle^{\tsigma^{\tilde{m}_\xi(\nu_\xi)_j}}
=(-1)^{\tilde{m}_\xi(\nu_\xi)_j}\}. \]
\begin{lemma} \label{unonpmlemma}
For $f\in\Theta_T^{F}$, $f\in\Theta_{T,\lambda}^{F}\Leftrightarrow
w_f\in X^{\underline{\nu}}_\inv$.
\end{lemma}
\begin{proof}
By the same argument as the proof of Lemma \ref{unonlemma}, we see that
\begin{equation*}
\epsilon_{T,f}(t)=
(-1)^{|\{(\xi,j)\in\Lambda_{w_f}^{1}(\underline{\nu})\,|\,
2\nmid\tilde{m}_\xi,2\nmid(\nu_\xi)_j,\alpha_{(\xi,j)}=-1\}|}.
\end{equation*} 
Hence $f\in\Theta_{T,\lambda}^{F}$ if and only if $w_f\in W(T)_\lambda^{F}$
and for all $(\xi,j)\in\Lambda_{w_f}^{1}(\underline{\nu})$,
\[ \tilde{d}_\xi=-1, 2\nmid(\nu_\xi)_j \Leftrightarrow
2\nmid\tilde{m}_\xi, 2\nmid(\nu_\xi)_j. \]
Clearly this is equivalent to $w_f\in X^{\underline{\nu}}_\inv$.
\end{proof}

\newpar
As in \S4.3,
\begin{equation*}
\begin{split}
\text{$\Fq$-rank}(T)+&\text{$\Fq$-rank}(Z_G((T\cap fKf^{-1})^{\circ}))\\
&\equiv \frac{n}{2}
+\ell_{w_f}^{1}(\underline{\nu})_0
+\ell_{w_f}^{2}(\underline{\nu})_{2\,\mathrm{mod}\,4}
+\frac{1}{2}\ell_{w_f}^{3}(\underline{\nu})_1\ \mathrm{mod}\ 2.
\end{split}
\end{equation*}
So Lusztig's formula gives (compare \S4.2):
\begin{equation*}
\begin{split}
\langle B_{\underline{\nu}}, \Ind_{K^{F}}^{G^{F}}(1)\rangle
&=\frac{1}{2}(-1)^{\frac{n}{2}}\sum_{w\in X^{\underline{\nu}}_\inv}
(-1)^{\ell_{w}^{1}(\underline{\nu})_0+
\ell_{w}^{2}(\underline{\nu})_{2\,\mathrm{mod}\,4}
+\frac{1}{2}\ell_{w}^{3}(\underline{\nu})_1}
2^{\ell_w^{1}(\underline{\nu})}\\
&\thickspace+\frac{1}{2}\epsilon\epsilon_{\underline{\nu}}
\sum_{w\in X^{\underline{\nu}}_{\ff-\inv}}
(-1)^{\ell_{w}^{2}(\underline{\nu})_{2\,\mathrm{mod}\,4}
+\frac{1}{2}\ell_{w}^{3}(\underline{\nu})_1}.
\end{split}
\end{equation*}

\newpar
Under the isomorphism $W(T)^{F}_\lambda\isomto\prod_\xi Z^{\nu_\xi}$ of \S1.4,
$X^{\underline{\nu}}_\inv$ corresponds to
\[ \{(w_\xi)\in\prod_{\xi\in\langle\tsigma\rangle\setminus L}
Z^{\nu_\xi}_\inv\,|\, \tilde{d}_\xi=-(-1)^{\tilde{m}_\xi}\Rightarrow
\ell_{w_\xi}^{1}(\nu_\xi)_1=0\}. \]
Hence
\begin{equation*}
\begin{split}
\langle B_{\underline{\nu}}, \Ind_{K^{F}}^{G^{F}}(1)\rangle
&=\frac{1}{2}(-1)^{\frac{n}{2}}\\
\times\prod_{\substack{\xi\in\langle\tsigma\rangle\setminus L\\
2\nmid\tilde{m}_\xi}}
&\left(\sum_{\substack{
w_\xi\in Z^{\nu_\xi}_\inv\\
\tilde{d}_\xi=1 \Rightarrow 
\ell_{w_\xi}^{1}(\nu_\xi)_1=0}}
\negthickspace\negthickspace\negthickspace\negthickspace\negthickspace
(-1)^{\ell_{w_\xi}^{1}(\nu_\xi)
+\ell_{w_\xi}^{2}(\nu_\xi)_{2\,\mathrm{mod}\,4}
+\frac{1}{2}\ell_{w_\xi}^{3}(\nu_\xi)_1}
2^{\ell_{w_\xi}^{1}(\nu_\xi)}\right)
\\
&\thickspace\times
\prod_{\substack{\xi\in\langle\tsigma\rangle\setminus L\\
2|\tilde{m}_\xi}}
\left(\sum_{\substack{
w_\xi\in Z^{\nu_\xi}_\inv\\
\tilde{d}_\xi=-1 \Rightarrow 
\ell_{w_\xi}^{1}(\nu_\xi)_1=0}}
\negthickspace\negthickspace
(-2)^{\ell_{w_\xi}^{1}(\nu_\xi)}\right)\\
&+\frac{1}{2}\epsilon
\prod_{\substack{\xi\in\langle\tsigma\rangle\setminus L\\
2\nmid\tilde{m}_\xi}}(-1)^{\frac{|\nu_\xi|}{2}}\epsilon_{\nu_\xi}
|Z^{\nu_\xi}_{\ff-\inv}|
\prod_{\substack{\xi\in\langle\tsigma\rangle\setminus L\\
2|\tilde{m}_\xi}}\epsilon_{\nu_\xi}|Z^{\nu_\xi}_{\ff-\inv}|.
\end{split}
\end{equation*}
In those factors of the second term for which $\tilde{m}_\xi$ is odd,
we have used the fact that if $w_\xi\in Z^{\nu_\xi}_{\ff-\inv}$, then
\[ \ell_{w_\xi}^{2}(\nu_\xi)_{2\,\mathrm{mod}\,4}+
\frac{1}{2}\ell_{w_\xi}^{3}(\nu_\xi)_1\equiv\frac{|\nu_\xi|}{2}
\ \mathrm{mod}\ 2. \]
So \eqref{unonpmeqn} follows by applying \eqref{combunoneqn},
\eqref{combunsoneqn}, \eqref{combglnoneqn}, \eqref{combglnsoneqn},
and \eqref{combglnspneqn}.
\section{Combinatorics of the Symmetric Group}
This section is devoted to the proof of the combinatorial facts
invoked in \S\!\S2-4. The notation introduced in \S1.2 will be used.
We say that a function $f$ on the set of partitions
is \emph{multiplicative} if
\[ f(\nu)=\prod_i f(i^{m_i(\nu)}),\ \forall \nu. \]
Examples of multiplicative functions of $\nu$ are $\epsilon_\nu$,
$z_\nu=|Z^{\nu}|$, and $|Z^{\nu}_\inv|$.

\newpar
Our first starting point is \cite[VII.(2.4)]{macdonald},
which as noted above is precisely
\begin{equation} \tag{\ref{combglnspneqn}}
|Z^{\nu}_{\ff-\inv}|
=\sum_{\substack{\rho\vdash|\nu|\\\rho\text{ even}}}
\chi_\nu^{\rho}.
\end{equation}
So to prove \eqref{combglnsoneqn},
it suffices to prove that
\begin{equation*}
\epsilon_\nu
\negthickspace\negthickspace
\sum_{\substack{w\in Z^{\nu}_\inv\\
\ell_w^{1}(\nu)_1=0}}
\negthickspace\negthickspace
(-2)^{\ell_{w}^{1}(\nu)}
=|Z^{\nu}_{\ff-\inv}|.
\end{equation*}
To see this, note that since both sides are multiplicative,
it suffices to consider the case
when $\nu$ is of the form $(a^{b})$, in which case both sides are
\[ \left\{\begin{array}{cl}
0,&\text{ if $a$ is odd and $b$ is odd,}\\
{\displaystyle a^{b/2}\frac{b!}{2^{b/2}(\frac{b}{2})!},}
&\text{ if $a$ is odd and $b$ is even, and}\\
{\displaystyle \sum_{r=0}^{\lfloor\frac{b}{2}\rfloor}
\binom{b}{2r}a^{r}\frac{(2r)!}{2^{r}r!},}
&\text{ if $a$ is even.}
\end{array}\right. \]
This fact also implies 
\begin{equation} \tag{\ref{combglnglneqn}}
\negthickspace\negthickspace
\sum_{(w,\epsilon)\in Z^{\nu}_{\star-\inv}}\negthickspace\negthickspace
\negthickspace\negthickspace
(-1)^{\ell_w^{2}(\nu)}
=\sum_{\substack{\rho\vdash|\nu|\\\text{$\rho$ even}}}
\chi_\nu^{\rho},
\end{equation}
once we note that
\begin{equation*}
\sum_{(w,\epsilon)\in Z^{\nu}_{\star-\inv}}\negthickspace\negthickspace
\negthickspace\negthickspace
(-1)^{\ell_w^{2}(\nu)}
=\sum_{\substack{w\in Z^{\nu}_\inv\\
\ell_w^{1}(\nu)_1=0}}
\negthickspace\negthickspace
(-1)^{\ell_w^{2}(\nu)}
2^{\ell_w^{1}(\nu)}
=\epsilon_\nu
\negthickspace\negthickspace
\sum_{\substack{w\in Z^{\nu}_\inv\\
\ell_w^{1}(\nu)_1=0}}
\negthickspace\negthickspace
(-2)^{\ell_w^{1}(\nu)},
\end{equation*}
since if there exists an involution $w\in Z^{\nu}$ fixing
no odd cycles of $\nu$, then $\ell(\nu)_1$ must be even, so
$(-1)^{\ell(\nu)}=\epsilon_\nu$.

\newpar
Our second starting point is \cite[I.8 Example 11]{macdonald}, which can
be rewritten:
\[ \sum_{w\in Z^{\nu}_\inv}
\negthickspace\negthickspace
(-1)^{\ell_w^{2}(\nu)}
=\sum_{\rho\vdash|\nu|}\chi_\nu^{\rho}. \]
So to prove \eqref{othercombglnoneqn},
it suffices to show that
\begin{equation*}
|\{w\in Z^{\nu}_\inv\,|\,\ell_w^{1}(\nu)_0=\ell_w^{2}(\nu)_0=0\}|
=\sum_{w\in Z^{\nu}_\inv}
\negthickspace\negthickspace
(-1)^{\ell_w^{2}(\nu)},
\end{equation*}
which we can again prove simply by observing that when $\nu=(a^{b})$ 
both sides are
\[ \left\{\begin{array}{cl}
0,&\text{ if $a$ is even and $b$ is odd,}\\
{\displaystyle a^{b/2}\frac{b!}{2^{b/2}(\frac{b}{2})!},}
&\text{ if $a$ is even and $b$ is even, and}\\
{\displaystyle \sum_{r=0}^{\lfloor\frac{b}{2}\rfloor}
\binom{b}{2r}a^{r}\frac{(2r)!}{2^{r}r!},}
&\text{ if $a$ is odd.}
\end{array}\right. \]

\newpar
The remaining identities require a different approach. We first prove
\begin{equation} \tag{\ref{combglnglnglneqn}}
\sum_{(w,\epsilon)\in Z^{\nu}_{(p^{+},p^{-})-\inv}}\negthickspace\negthickspace
(-1)^{\ell_w^{2}(\nu)}
=\sum_{\rho\vdash|\nu|} |\T_{(p^{+},p^{-})}(\rho')|\chi_\nu^{\rho}.
\end{equation}
By definition of induction product
(\cite[I.7]{macdonald}), and using
\eqref{combglnspneqn},
\begin{equation*}
\begin{split}
\sum_{(w,\epsilon)\in Z^{\nu}_{(p^{+},p^{-})-\inv}}&
\negthickspace\negthickspace
(-1)^{\ell_w^{2}(\nu)}\\
&=\sum_{r=\lceil\frac{p^{+}+p^{-}-|\nu|}{2}\rceil}
^{\mathrm{min}\{p^{+},p^{-}\}}
\left((\sum_{\substack{
\mu\vdash|\nu|+2r-p^{+}-p^{-}\\ \mu'\text{ even}}} \chi^{\mu})
.\chi^{(p^{+}-r)}.\chi^{(p^{-}-r)}\right)(w_\nu) \\
&=\sum_{r=\lceil\frac{p^{+}+p^{-}-|\nu|}{2}\rceil}
^{\mathrm{min}\{p^{+},p^{-}\}}
\sum_{\rho\vdash|\nu|} b(p^{+}-r,p^{-}-r,\rho')\chi_\nu^{\rho},
\end{split}
\end{equation*}
where, by Pieri's formula (\cite[(5.16)]{macdonald}),
$b(p^{+}-r,p^{-}-r,\rho')$ is the number of ways of removing first
a vertical $(p^{-}-r)$-strip, then a vertical $(p^{+}-r)$-strip, from
the Young diagram of $\rho'$, to leave a diagram with all rows of
even length. Now every signed tableaux $T\in\T_{(p^{+},p^{-})}(\rho')$
determines uniquely an $r$ as above
and such a way of removing strips, as follows:
\begin{itemize}
\item  order rows of equal length so that rows ending
$\boxminus$ are below rows ending $\boxplus$;
\item  take the vertical $(p^{-}-r)$-strip to consist of all final
boxes signed $\boxminus$;
\item  take the vertical $(p^{+}-r)$-strip to consist of all final
boxes signed $\boxplus$ in rows of odd length, including those made odd
by removal of the first strip.
\end{itemize}
This correspondence is clearly bijective, which proves
\eqref{combglnglnglneqn}.

\newpar
The proof of
\begin{equation} \tag{\ref{combunununeqn}}
\sum_{(w,\epsilon)\in Z^{\nu}_{(p^{+},p^{-})-\inv}}
\negthickspace\negthickspace
(-1)^{\ell_w^{2}(\nu)_{0\,\mathrm{mod}\,4}+\frac{1}{2}\ell_w^{3}(\nu)_1}
=\sum_{\rho\vdash|\nu|} (-1)^{n(\rho)}
|\T_{(p^{+},p^{-})}(\rho')^{\psi}|\chi_\nu^{\rho}
\end{equation}
proceeds similarly:
\begin{equation*}
\begin{split}
&\sum_{(w,\epsilon)\in Z^{\nu}_{(p^{+},p^{-})-\inv}}
\negthickspace\negthickspace
(-1)^{\ell_w^{2}(\nu)_{0\,\mathrm{mod}\,4}+\frac{1}{2}\ell_w^{3}(\nu)_1}\\
&=\sum_{r=\lceil\frac{p^{+}+p^{-}-|\nu|}{2}\rceil}
^{\mathrm{min}\{p^{+},p^{-}\}}
\left(((-1)^{\frac{|\nu|+2r-p^{+}-p^{-}}{2}}\sum_{\substack{
\mu\vdash|\nu|+2r-p^{+}-p^{-}\\ \mu'\text{ even}}} \chi^{\mu})
.\chi^{(p^{+}-r)}.\chi^{(p^{-}-r)}\right)(w_\nu) \\
&=\sum_{\rho\vdash|\nu|}
c(p^{+},p^{-},\rho')\chi_\nu^{\rho},
\end{split}
\end{equation*}
where by the same bijection as before,
\begin{equation*}
\begin{split}
c&(p^{+},p^{-},\rho')\\
&=\sum_{T\in\T_{(p^{+},p^{-})}(\rho')}
\negthickspace\negthickspace
(-1)^{\begin{array}{l}
{\scriptstyle\frac{1}{2}(|\rho|-
|\{\text{rows of $T$ ending $\boxminus$}\}|}\\
{\scriptstyle\ -|\{\text{odd rows of $T$ ending $\boxplus$}\}|
-|\{\text{even rows of $T$ ending $\boxminus$}\}|)}
\end{array}}\\
&=(-1)^{\frac{1}{2}(|\rho|-\ell(\rho')_1)}
\negthickspace\negthickspace
\sum_{T\in\T_{(p^{+},p^{-})}(\rho')}
\negthickspace\negthickspace
(-1)^{|\{\text{even rows of $T$ ending $\boxminus$}\}|}.
\end{split}
\end{equation*}
Now
\begin{equation*}
(-1)^{\frac{1}{2}(|\rho|-\ell(\rho')_1)}=
(-1)^{\ell(\rho')_{2\,\mathrm{mod}\,4}
+\ell(\rho')_{3\,\mathrm{mod}\,4}}
= (-1)^{\sum_i \binom{\rho_i'}{2}} = (-1)^{n(\rho)},
\end{equation*}
and by grouping together signed tableaux which differ only in even rows
it is easy to see that the sum equals $|\T_{(p^{+},p^{-})}(\rho')^{\psi}|$.
So \eqref{combunununeqn} is proved.

\newpar
Our next task is to modify this proof of \eqref{combunununeqn} to derive
\begin{equation} \tag{\ref{combununeqn}}
\begin{split}
\sum_{(w,\epsilon)\in Z^{\nu}_{\star-\inv}}\negthickspace\negthickspace
(-1)&^{l_w^{2}(\nu)_{0\,\mathrm{mod}\,4}+\frac{1}{2}\ell_w^{3}(\nu)_1}\\
&=\sum_{\substack{\rho\vdash|\nu|\\2|m_{2i+1}(\rho')}} 
(-1)^{n(\rho)}
(\prod_i (m_{2i}(\rho')+1))\chi_\nu^{\rho}.
\end{split}
\end{equation}
We may assume that $|\nu|$ is even, for otherwise both sides vanish.
Following the above pattern,
we need to replace $\chi^{(p^{+}-r)}.\chi^{(p^{-}-r)}$ with the class
function on $S_{|\nu|-2r}$ defined by
\[ w\mapsto|\{(A^{+},A^{-})\,|\,\{1,\cdots,|\nu|-2r\}=A^{+}
{\textstyle\coprod} A^{-},
|A^{+}|=|A^{-}|, w(A^{\pm})=A^{\mp}\}|. \]
It is easy to see that this is
$\sum_{i=0}^{\frac{|\nu|}{2}-r} (-1)^{i}\chi^{(|\nu|-2r-i,i)}$.
Thus
\begin{equation*}
\begin{split}
&\sum_{(w,\epsilon)\in Z^{\nu}_{\star-\inv}}
\negthickspace\negthickspace
(-1)^{\ell_w^{2}(\nu)_{0\,\mathrm{mod}\,4}+\frac{1}{2}\ell_w^{3}(\nu)_1}\\
&=\sum_{r=0}
^{\frac{|\nu|}{2}}
\left(((-1)^{r}\sum_{\substack{
\mu\vdash 2r\\ \mu'\text{ even}}} \chi^{\mu})
.(\sum_{i=0}^{\frac{|\nu|}{2}-r} (-1)^{i}\chi^{(|\nu|-2r-i,i)})\right)
(w_\nu) \\
&=\sum_{\rho\vdash|\nu|}
d(\rho')\chi_\nu^{\rho},
\end{split}
\end{equation*}
where, analogously to the above, we can write $d(\rho')$ as the sum,
over $T\in\T_{(0,0)}(\rho')$, of a sign determined by $T$. As in the proof
of \eqref{combunununeqn}, the $(-1)^{r}$ contribution to the sign is
\[ (-1)^{n(\rho)+|\{\text{even rows of $T$ ending $\boxminus$}\}|}. \]
The $(-1)^{i}$ contribution is trickier to rephrase in terms of $T$, but
an examination of the bijection
\[ \mathrm{Tab}(\rho'-\mu',(\frac{|\nu|}{2},\frac{|\nu|}{2}))
\isomto \coprod_{i=0}^{\frac{|\nu|}{2}-r}
\mathrm{Tab}^{0}(\rho'-\mu',(|\nu|-2r-i,i)) \]
defined in \cite[(9.4)]{macdonald} (in the proof of the Littlewood-Richardson
Rule) reveals that the correct reformulation is
\[ (-1)^{|\{\text{rows of $T$ ending $\boxminus$}\}|-m(T)}, \]
where $m(T)$ is the maximum, over all rows $R$ of $T$, of the quantity
\begin{equation*}
\begin{split}
&|\{\text{odd rows ending $\boxminus$ below or equal to $R$}\}|\\
&\qquad\qquad\qquad\qquad\qquad\qquad
-|\{\text{odd rows ending $\boxplus$ below or equal to $R$}\}|
\end{split}
\end{equation*}
(assuming that the rows of $T$ are ordered, as above, so that rows
ending $\boxminus$ come below rows ending $\boxplus$ of the same length).
Hence
\[ d(\rho')=(-1)^{n(\rho)}\negthickspace
\sum_{T\in\T_{(0,0)}(\rho')}
\negthickspace\negthickspace
(-1)^{|\{\text{odd rows of $T$ ending $\boxminus$}\}|-m(T)}, \]
and we are reduced to proving that
\begin{equation*}
\begin{split}
\sum_{T\in\T_{(0,0)}(\rho')}&
\negthickspace\negthickspace
(-1)^{|\{\text{odd rows of $T$ ending $\boxminus$}\}|-m(T)}\\
&=\left\{\begin{array}{cl}
{\displaystyle \prod_i(m_{2i}(\rho')+1),}
&\text{ if $2\,|\,m_{2i+1}(\rho')$, $\forall i$}\\
0, &\text{ otherwise.}
\end{array}\right.
\end{split}
\end{equation*}
Grouping together signed tableaux which differ only in even rows,
we see that we may assume that $\rho'$ has only odd parts. In lieu
of a direct proof, we can deduce this from \eqref{combglnglneqn}, proved above.
It says that
\begin{equation*}
\begin{split}
\sum_{\substack{\rho\vdash|\nu|\\\rho\text{ even}}}\chi_\nu^{\rho}
&= \sum_{(w,\epsilon)\in Z^{\nu}_{\star-\inv}}\negthickspace\negthickspace
(-1)^{\ell_w^{2}(\nu)}\\
&= \sum_{r=0}
^{\frac{|\nu|}{2}}
\left((\sum_{\substack{
\mu\vdash 2r\\ \mu'\text{ even}}} \chi^{\mu})
.(\sum_{i=0}^{\frac{|\nu|}{2}-r} (-1)^{i}\chi^{(|\nu|-2r-i,i)})\right)(w_\nu),
\end{split}
\end{equation*}
whence
\[ \sum_{T\in\T_{(0,0)}(\rho')}
\negthickspace\negthickspace
(-1)^{|\{\text{rows of $T$ ending $\boxminus$}\}|-m(T)}
=\left\{\begin{array}{cl}
1, &\text{ if $2\,|\,m_{i}(\rho')$, $\forall i$}\\
0, &\text{ otherwise.}
\end{array}\right. \]
When $\rho'$ has only odd parts, this is precisely the statement we want.

\newpar
It is now easy to prove \eqref{combunsoneqn}, since
\begin{equation*}
\begin{split}
\sum_{\substack{w\in Z^{\nu}_\inv\\
\ell_w^{1}(\nu)_1=0}}
&\negthickspace\negthickspace
(-1)^{\ell_w^{1}(\nu)+\ell_w^{2}(\nu)_{2\,\mathrm{mod}\,4}
+\frac{1}{2}\ell_w^{3}(\nu)_1}
2^{\ell_w^{1}(\nu)}\\
&=(-1)^{\ell(\nu)}\negthickspace\negthickspace
\sum_{(w,\epsilon)\in Z^{\nu}_{\star-\inv}}
\negthickspace\negthickspace
(-1)^{\ell_w^{2}(\nu)_{0\,\mathrm{mod}\,4}
+\frac{1}{2}\ell_w^{3}(\nu)_1}\\
&=\epsilon_\nu\negthickspace\negthickspace
\sum_{\substack{\rho\vdash|\nu|\\2|m_{2i+1}(\rho')}}\negthickspace
(-1)^{n(\rho)}
(\prod_i (m_{2i}(\rho')+1))
\chi_\nu^{\rho}\\
&=\sum_{\substack{\rho\vdash|\nu|\\2|m_{2i+1}(\rho)}}\negthickspace
(-1)^{n(\rho')}
(\prod_i (m_{2i}(\rho)+1))
\chi_\nu^{\rho},
\end{split}
\end{equation*}
as required. Finally, \eqref{combglnoneqn}
follows from \eqref{combglnglnglneqn} by summing
over all signatures:
\begin{equation*}
\begin{split}
\sum_{w\in Z^{\nu}_\inv}\negthickspace\negthickspace
(-2)^{\ell_w^{1}(\nu)}
&=(-1)^{\ell(\nu)}\negthickspace\negthickspace
\sum_{(w,\epsilon)\in Z^{\nu}_{\pm-\inv}}
\negthickspace\negthickspace
(-1)^{\ell_w^{2}(\nu)}\\
&=(-1)^{|\nu|}\epsilon_\nu\sum_{\rho\vdash|\nu|}
|\T_{\pm}(\rho')|\chi_\nu^{\rho}\\
&=(-1)^{|\nu|}\sum_{\rho\vdash|\nu|}
(\prod_i (m_i(\rho)+1))\chi_\nu^{\rho},
\end{split}
\end{equation*}
and \eqref{combunoneqn} follows from \eqref{combunununeqn} in an
analogous way:
\begin{equation*}
\begin{split}
\sum_{w\in Z^{\nu}_\inv}
\negthickspace\negthickspace
(-1)^{\ell_w^{1}(\nu)_0+\ell_w^{2}(\nu)_{2\,\mathrm{mod}\,4}
+\frac{1}{2}\ell_w^{3}(\nu)_1}
&2^{\ell_w^{1}(\nu)}
=\epsilon_\nu\sum_{\rho\vdash|\nu|}(-1)^{n(\rho)}
|\T_{\pm}(\rho')^{\psi}|\chi_\nu^{\rho}\\
&=\sum_{\substack{\rho\vdash|\nu|\\2|m_{2i}(\rho)}}
(-1)^{n(\rho')}(\prod_i(m_{2i+1}(\rho)+1))\chi_\nu^{\rho}.
\end{split}
\end{equation*}

\end{document}